\documentclass[10pt,a4paper]{article}
\usepackage[utf8]{inputenc}
\usepackage[margin=2cm]{geometry}
\usepackage{bm} % for bold maths symbols
\usepackage{enumitem} % for more flexibility with lists
\usepackage{array}  % to give vertical centering in tables
\usepackage[pdftex]{graphicx} % for including images
\usepackage{subcaption} % for side by side figures
\usepackage{amssymb} % for mathbb
\usepackage{amsmath} % for align and tfrac
\usepackage{amsthm} % for theorem / proofs
\usepackage{natbib}
\usepackage{hyperref} % for interactable links
\usepackage{authblk} % for author and affiliations
\usepackage{fancyhdr} % For headers and footers
\usepackage{tikz-cd}
\newcommand{\grad}[1]{{\bm{\nabla}{#1}}}

\newcommand{\pfrac}[2]{\frac{\partial{#1}}{\partial{#2}}}
\newcommand{\dx}[1]{\hspace{0.75mm}\mathrm{d}{#1}}

\newcommand{\curlybrackets}[1]{\left\lbrace{#1}\right\rbrace}

\newtheorem{requirement}{Requirement}

\fancypagestyle{plain}{
\fancyhf{}
\lfoot{\textcopyright \ Crown Copyright, Met Office}
\cfoot{\thepage}
}

\title{Physics-Dynamics-Chemistry Coupling Across Different Meshes in LFRic-Atmosphere: Formulation and Idealised Tests}

\author[1,2]{Alex Brown}
\author[1]{Thomas M. Bendall\thanks{Corresponding author, contact at thomas.bendall@metoffice.gov.uk}}
\author[1]{Ian Boutle}
\author[1]{Thomas Melvin}
\author[1]{Ben Shipway}
\affil[1]{Dynamics Research, Met Office, Exeter, UK}
\affil[2]{Department of Mathematics, University of Exeter, UK}
\date{}

\begin{document}

\maketitle

\begin{abstract}
\noindent The main components of an atmospheric model for numerical weather prediction are the \textit{dynamical core}, which describes the resolved flow, and the \textit{physical parametrisations}, which capture the effects of unresolved processes.
Additionally, models used for air quality or climate applications may include a component that represents the evolution of chemicals and aerosols within the atmosphere.
While traditionally all these components use the same mesh with the same resolution,
we present a formulation for the different components to use a series of nested meshes, with different horizontal resolutions.
This gives the model greater flexibility in the allocation of computational resources,
so that resolution can be targeted to those parts which provide the greatest benefits in accuracy.
\\
\\
The formulation presented here concerns the methods for mapping fields between meshes, and is designed for the compatible finite element discretisation used by LFRic-Atmosphere, the Met Office's next-generation atmosphere model.
Key properties of the formulation include the consistent and conservative transport of tracers on a mesh that is coarser than the dynamical core, and the handling of moisture to ensure mass conservation without generation of unphysical negative values.
Having presented the formulation, it is then demonstrated through a series of idealised test cases which show the feasibility of this approach.
\end{abstract}

\section{Introduction} \label{sec:intro}
Due to the complexity of the equations that describe the evolution of the atmosphere, the numerical models typically used in simulating the weather and climate are broken down into different components, each describing different processes.
The \textit{dynamical core} (or ``dynamics'') discretises the equations for resolved fluid motions.
The \textit{physical parametrisations} (or ``physics'') capture the non-fluid processes and the non-resolved fluid processes.
A final component, most often found in climate and air quality models, describes the transport of aerosols and chemicals, and the reactions between the chemicals (this component will be referred to collectively as ``chemistry'' throughout). 
As discussed by \citet{gross2018physics}, these components are generally written independently from one another, but coupled together in some way to form the whole atmospheric model.
This structure has often evolved naturally, as the complexity of the equations governing the Earth system necessitates different terms being discretised and evaluated separately. \\
\\
Traditionally, in numerical weather prediction (NWP) and climate models the dynamical core, physical parametrisations and the chemistry component are computed on the same mesh, and often this choice has been made to simplify the coupling between the different components.
Notable exceptions to this include those models which use spectral element or spectral transform methods in their dynamical core, such as ECMWF's IFS model (\citep{roberts2018climate},\citet{malardel2016new}); NCAR's CAM-SE spectral element model, in which \cite{herrington2019exploring} and \citet{herrington2019physics} have recently explored the use of a coarser physics grid; and the Department of Energy's E3SM spectral element model, in which \cite{hannah2021separating} and
\cite{bradley2022islet} have investigated the use of alternative physics and tracer transport grids.
Another related endeavour is a climate configuration of the Met Office's Unified Model (UM), known as Junior-Senior, which is motivated by reducing the large computational cost of the chemistry component \citep{stringer2018hybrid}.
\\
\\
This work explores removing the assumption that the different atmospheric components use the same grid, in the context of the Met Office's next-generation LFRic-Atmosphere model, which uses a compatible finite element discretisation in its dynamical core, GungHo.
In this paper we present a formulation for coupling together the different resolution dynamics, physics and chemistry components, inspired by the approach of \citet{herrington2019exploring} and \citet{herrington2019physics}.
The formulation is then tested through a series of idealised examples.
Future work will seek to investigate and understand the consequences of this new capability within full NWP and climate models.

\subsection{Background and Motivation} \label{sec:background}

Whilst traditionally the components of atmospheric models use a mesh of the same resolution, the concept of using different meshes for different components is not novel. However, there are contrasting arguments for how the resolution of the physical parametrisations should be changed relative to that of the dynamical core.
\\
\\
\cite{gross2018physics} presents both arguments. On the one hand, it is argued that computing physics on a high-resolution mesh means sampling the fields from the dynamical core more finely, comparing this to the ``subcolumns" approach that is used in some cloud-schemes. The physical parametrisations generally describe non-linear processes, so computing these at a higher resolution may give a better representation of their effect on the resolved flow. One example of a model that uses a higher resolution for the physical parametrisations is ECMWF's spectral IFS model \citep{roberts2018climate}. The grid on which physical parametrisations are computed on has more degrees of freedom than the number of wave modes used in the spectral part of the dynamical core, which  \cite{malardel2016new} found to give reduced aliasing and better mass conservation. Whether the same benefits would apply to non-spectral models is not clear.
\\
\\
On the other hand, as argued by \cite{lander1997believable}, the physics should perhaps only be passed well-resolved ``believable" scales from the dynamics, as the numerical errors in the solutions may be amplified by the non-linear physical parametrisations. These numerical errors are likely to be largest at the smallest scales of motion, which are generally poorly-resolved by the dynamical core. Therefore, by computing the physics at coarser resolution to dynamics, the physical parametrisations only act upon fields from which these poorly resolved scales have been filtered. This approach of \cite{lander1997believable} was also considered in the context of a spectral transform model.
\\
\\
A final factor relates to that of computational cost. If the different components of the atmospheric model can use different grids, then the computational resources can be targeted to the part of the model that provides the greatest benefit. Alternatively, there may be parts of the model whose resolution can be reduced without particularly degrading the solution quality, freeing up computational resources to be assigned elsewhere, possibly into increasing model complexity rather than resolution.
\\
\\
Some of these ideas have been explored in the spectral element CAM-SE model. \cite{herrington2019physics} implemented an alternative quasi-equal-area finite volume physics grid, which had the same number of degrees of freedom as the dynamics grid, and on which tracer advection is also computed. This reduced grid imprinting and spurious vertical velocity noise over orography. \cite{herrington2019exploring} extended this to use a coarser physics grid, with a $5/9$ reduction in the number of columns in the physics grid with respect to the dynamics grid, while the tracer grid remained at the same effective resolution as the dynamics grid. In \citet{herrington2019exploring} and \citet{herringon2019physics}, the model's prognostic variables were mapped from the dynamics grid to the other grids, while only increments were mapped back to the dynamics grid. Momentum components were interpolated by evaluating their basis functions at the physics degrees of freedom, while pressure and temperature variables were integrated over the coarse control volumes. Tracers are mapped to the physics grid by a high-order reconstruction that preserves tracer shape, linear correlations and conserves mass. Increments were mapped with an alternative algorithm which alters the mixing ratio increment in order to also preserve shape, linear correlations, conserve mass, as well as maintaining consistency and positivity. \cite{herrington2019exploring} demonstrated that the effective resolution was not degraded through aquaplanet simulations, allowing for future computational savings. It also reduced noise over steep orography at element boundaries, a common problem in spectral element models, shown through a Held-Suarez test with real orography.
\\
\\
\cite{hannah2021separating} expanded on this approach in the spectral element E3SM model, investigating the use of a higher-resolution mesh for the physics parametrisations, but found no qualitative benefit. When the physics parametrisation mesh was lower-resolution, \cite{hannah2021separating} showed no degradation in the solution for a simulated climate. The lower-resolution physics grid was further shown to reduce grid imprinting over orography and demonstrated significant computational savings. \cite{bradley2022islet} extended this by implementing an alternative grid for tracer transport, using a interpolation semi-Lagrangian finite element transport scheme, in the E3SM model where physics and chemistry are on the lower-resolution physics grid.
\\
\\
The Met Office's hybrid-resolution version of the UKESM earth system model \citep{stringer2018hybrid} (Junior-Senior) runs a high-resolution version of the Met Office's Unified Model (UM) without UKCA chemistry and aerosol (dynamics and physics) driving a low-resolution version of the UM with UKCA (dynamics, physics and chemistry). This is motivated by reducing the significant computational cost of the chemistry model, which contains a significant number of chemical and aerosol species. The formulation presented in this paper could be used to address the same problem in the new LFRic-Atmosphere model.
\\
\\
The work presented here has taken significant inspiration from \cite{herrington2019physics} and \cite{herrington2019exploring}, but differs in some key elements. Whereas those works used a spectral element dynamical core, this work considers a dynamical core with the lowest-order compatible finite element discretisation of \cite{melvin2019mixed}. This has implications for the staggering of the different prognostic variables and hence the operators used to map these fields. In this work, the mesh used for the physical parametrisations has the same structure as the dynamical core (but of different resolution and with cells exactly nested within or exactly nesting those of the dynamical core), whereas the approach in CAM-SE used a finite volume grid overlaying the spectral element grid. The two approaches also preserve similar but subtly different properties.

\subsection{Scope} \label{sec:scope}
The formulation presented in this work is designed for a model with three constituent parts: a dynamical core, coupled to a chemistry and aerosol model, and a set of physical parametrisations.
The dynamical core evolves a set of dynamical prognostic variables (including fields describing the moist composition of the atmosphere), while the chemistry and aerosol component evolves a different set of variables, describing the chemical and aerosol species in the atmosphere.
The physical parametrisations provide updates to the dynamical prognostic variables and the chemical and aerosol species, but also depend on a set of prescribed auxiliary variables.
The chemicals and aerosols do not feed directly back into the dynamical core, but they may appear as auxiliary fields to the physics schemes.
\\
\\
Following the motivations laid out earlier in Section \ref{sec:background}, the formulation is designed for three different types of interaction between these components:
\begin{enumerate}[leftmargin=*]
\item physical parametrisations that are computed on a finer mesh than the dynamical core;
\item physical parametrisations that are computed on a coarser mesh than the dynamical core;
\item a chemistry and aerosol component (including tracer transport) computed on a coarser mesh than the dynamical core.
\end{enumerate}
The interactions between components that use different meshes involve mapping fields from one mesh to another.
To avoid complications relating to the averaging of vector-valued fields, only physics parametrisations providing updates to scalar-valued fields are computed on a different mesh to the dynamical core.
\\
\\
The choices of mesh for these components are constrained by some crucial simplifications.
The three-dimensional meshes are extruded, so that they are the product of a two-dimensional horizontal mesh with a vertical one-dimensional mesh.
The two-dimensional horizontal mesh consists of quadrilateral cells, resulting in hexahedral cells in the three-dimensional mesh.
All the components use meshes with the same vertical structure, so that the resolution only differs in the horizontal part.
Cells on a finer mesh are exactly nested within those of a coarser mesh, which offers two significant design advantages.
Firstly, it is straightforward to calculate the size of the overlapping region between cells on two different meshes.
Secondly, this facilitates an efficient parallel distribution of memory
so that data corresponding to fields on different meshes can be geographically distributed in the same way, minimising the amount of data communication required to map fields from one mesh to another.
\\
\\
The purpose of this paper is to present the formulation used in LFRic-Atmosphere for coupling together these components when they use meshes of different resolutions.
The approach is demonstrated through a series of idealised test cases, which illustrate various aspects of the formulation.
In particular we focus on the transport of tracers on a coarser mesh, and at this stage do not demonstrate a dynamical core coupled to a full suite of physical parametrisations or a chemistry and aerosol model.  
Future work will extend this approach to full NWP and climate configurations and will explore the consequences of different choices of mesh for individual physics schemes and the subsequent consequences on the model's performance.
\\
\\
The remainder of the paper is organised as follows.
Section \ref{sec:prelims} specifies the prognostic variables used by the model, and also sets out the notation used in this paper to describe the formulation for coupling components of different resolutions.
Then, Section \ref{sec:desirable properties} discusses the properties of the formulation that we consider to be important.
The formulation, including the specific operators for mapping fields between meshes, is presented in Section \ref{sec:formulation}, which also shows that these operators satisfy the properties of Section \ref{sec:desirable properties}.
Section \ref{sec:results} demonstrates the formulation through idealised test cases.

\section{Preliminaries} \label{sec:prelims}

\subsection{Prognostic Variables} \label{sec:prognostics}
LFRic-Atmosphere's dynamical core, called GungHo, solves for the wind velocity $\bm{u}$, the dry density $\rho_d$, the Exner pressure $\Pi$ and the (dry) potential temperature $\theta$.
There are $N_r$ species of moisture which are described through mass mixing ratios, with the $r$-th species given by $m_r:=\rho_r/\rho_d$, where $\rho_r$ is a moisture density.
Collectively these prognostic variables can be described as a single state vector $\bm{X}$,
\begin{equation}
\bm{X} = \left(\bm{u}, \rho_d, \Pi, \theta, m_1, \dots, m_{N_r}\right).
\end{equation}
The mass mixing ratio $a_Y:=\rho_Y/\rho_d$ is also used to represent the $Y$-th chemical/aerosol species, so that if the model evolves $N_Y$ species, then the vector $\bm{Y}$ can be used for the chemical and aerosol species:
\begin{equation}
\bm{Y} = \left(a_1,\dots, a_{N_Y}\right).
\end{equation}
The model solves the compressible Euler equations, with additional equations for the moisture, chemical and aerosol variables:
\begin{subequations} \label{eqn:compressible_euler}
\begin{align}
& \pfrac{\bm{u}}{t}+\left(\bm{u\cdot\nabla}\right)\bm{u} + 2\bm{\varOmega}\times\bm{u} + \frac{c_p\theta(1 + m_v R_v/R_d)}{1+\sum_{r=1}^{N_r} m_r}\grad{\Pi}+\bm{g} = \bm{S}_u, \\
& \pfrac{\rho_d}{t} + \bm{\nabla\cdot}\left(\rho_d\bm{u}\right) = 0 \label{eqn:compressible_euler_rho}, \\
& \pfrac{\theta}{t} + \left(\bm{u\cdot\nabla}\right)\theta = S_\theta, \\
& \pfrac{m_r}{t} + \left(\bm{u\cdot\nabla}\right)m_r = S_r, \quad r\in[1,N_r], \\
& \pfrac{a_Y}{t} + \left(\bm{u\cdot\nabla}\right)a_Y = S_Y, \quad Y\in[1,N_Y]  \label{eqn:compressible_euler_tracer},
\end{align}
\end{subequations}
where $\bm{S}_u$, $S_\theta$ and $S_r$ represent the changes to the prognostic variables that are computed through the physical parametrisations.
The $S_Y$ variables describes sources, sinks and reactive effects computed by the chemical and aerosol model.
Equation \eqref{eqn:compressible_euler} is supplemented by the equation of state for an ideal gas,
\begin{equation}
\Pi = \left(\frac{\rho_d R_d \theta (1 + m_v R_v/R_d)}{p_0}\right)^\frac{R_d}{c_p+R_d},
\end{equation}
with $m_v$ as the mixing ratio of water vapour.
The constants are: the specific gas constant for dry air $R_d$, the specific gas constant for water vapour $R_v$, the specific heat capacity of dry air at constant pressure $c_p$, the reference pressure $p_0$, the gravitational field vector $\bm{g}$ and the Earth's rotation vector $\bm{\varOmega}$.
\subsection{Overview of LFRic-Atmosphere} \label{sec:lfric}
LFRic-Atmosphere is the Met Office's new weather forecasting and climate model, designed to exploit the next generation of supercomputers, as described in \cite{adams2019lfric}.
A major issue for adapting the Met Office's Unified Model (UM) \citep{wood2014inherently,walters2017met} to these supercomputers is the latitude-longitude mesh used for global simulations by the UM's dynamical core, ENDGame.
The latitude-longitude mesh has a convergence of spatial points at the poles, which leads to a bottleneck in data communication and a resolution gap between the poles and the equator.
This presents an unsustainable constraint on the UM's scalability as horizontal resolution is increased.
\\
\\
Key to LFRic's design is the use of a quasi-uniform cubed-sphere mesh, in both the physical parametrisations and the dynamical core, GungHo. ENDGame used C-grid and Charney-Phillips staggerings to obtain good linear wave dispersion properties and to avoid computational modes It was been shown by  \cite{cotter2012mixed}, \cite{cotter2014finite} and  \cite{thuburn2015primal}  that a compatible finite element discretisation can replicate these desirable properties, while also facilitating the move to a non-orthogonal mesh.
. %! This paragraph is a little clunky, I think we could make it flow better.
\\
\\
In the compatible finite element discretisation used by GungHo, all of the prognostic variables are discretised as a sum of coefficients multiplying basis functions, with the basis functions localised to a single element or set of elements surrounding a cell edge or vertex.
A \textit{finite element} is described by the choice of basis functions (usually polynomials) and their continuity between cells; then the combination of a finite element with the model's mesh defines the \textit{function space}.
In a \textit{compatible} finite element discretisation, variables lie in function spaces that form a de Rham complex, so that the vector calculus relationships between the discretised variables mimic those from the continuous equations.
A formal discussion of these concepts can be found in \cite{arnold2010finite} and \citet{cotter2023compatible}. \\
\\
GungHo uses the lowest-order finite elements of the Raviart-Thomas de Rham complex, that are extended to hexahedral cells through a tensor-product construction.
In this compatible finite element set-up, the prognostic variables are contained within three function spaces: $\mathbb{V}_u$, $\mathbb{V}_\theta$ and $\mathbb{V}_\rho$ (with the subscript denoting the variables contained within those spaces). 
The DoFs of $\mathbb{V}_\rho$ lie at the centre of cells, which corresponds to basis functions that are constant within a cell (and discontinuous between cells).
The Arakawa C-grid is replicated by staggering the DoFs of $\mathbb{V}_u$ from those of $\mathbb{V}_\rho$, so that the the DoFs of $\mathbb{V}_u$ are located at the faces of cells.
Then the values of fields at the $\mathbb{V}_u$ DoFs represent the normal fluxes of that field through the faces of the element.
The compatibility of $\mathbb{V}_u$ and $\mathbb{V}_\rho$ means that for any $\bm{u}\in\mathbb{V}_u$, then $\bm{\nabla\cdot u}\in\mathbb{V}_\rho$.
The DoFs of $\mathbb{V}_\theta$ are co-located with the vertical component of $\mathbb{V}_u$, and so the DoFs are located at the centre of the top or bottom surfaces of cells, which was shown by \citet{melvin2018choice} to mimic the Charney-Phillips staggering.
More description of these spaces is given by \citet{melvin2019mixed} and \citet{bendall2020compatible}, while representations of them are displayed in Table \ref{tab:elements}. \\
\\
With these function spaces, \eqref{eqn:compressible_euler} is discretised by taking $\bm{u}\in\mathbb{V}_u$ and $\rho_d,\Pi\in\mathbb{V}_\rho$.
In this work we consider chemicals and aerosol variables with mixing ratios $a_Y\in\mathbb{V}_\rho$, although the formulation in Section \ref{sec:formulation} can be extended to the case of $a_Y\in\mathbb{V}_\theta$.
The moisture variables are co-located with $\theta$, so that $\theta,m_r\in\mathbb{V}_\theta$, to give an accurate representation of the saturation curve and the latent heat exchanges associated with changes of phase.
\begin{table}[h!]
\centering
\begin{tabular}{c|c|c|c|c}
Space       & $\mathbb{V}_u$ & $\mathbb{V}_\theta$ & $\mathbb{V}_\rho$ & $\widetilde{\mathbb{V}}_\rho$
\\
\hline
Variables   & $\bm{u}$ & $\theta$, $m_r$ & $\rho_d$, $\Pi$, $a_Y$ & $\widetilde{\rho}_r$ \\
\hline
            &  \includegraphics[width=0.18\textwidth]{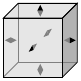} &
               \includegraphics[width=0.18\textwidth]{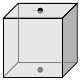} &
               \includegraphics[width=0.18\textwidth]{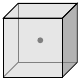} &
               \includegraphics[width=0.18\textwidth]{tab_1_fig_c_W3.pdf}

\end{tabular}
\caption{The finite elements used by GungHo in the discretisation of its prognostic variables.
The spaces $\mathbb{V}_u$ and $\mathbb{V}_\rho$ form part of a de Rham complex, so that if $\bm{u}\in\mathbb{V}_u$ then $\bm{\nabla\cdot u}\in\mathbb{V}_\rho$.
The degrees of freedom for $\mathbb{V}_u$ correspond to the fluxes through each face of the hexahedron, while there is one degree of freedom per cell for $\mathbb{V}_\rho$, representing the field's value at the cell's centre.
The degrees of freedom for $\mathbb{V}_\theta$ are in the centre of the top and bottom faces of cells.
The density of the $r$-th moisture species, $\widetilde{\rho}_r$ is described using the same elements as $\mathbb{V}_\rho$ but on a vertically-shifted mesh.} \label{tab:elements}
\end{table}

\subsection{Moisture conservation} \label{sec:moisture_conservation}
With $m_r\in\mathbb{V}_\theta$, conservation of the mass of moisture requires more steps than if it were located in $\mathbb{V}_\rho$.
As described by \cite{bendall2023solution}, this is addressed in GungHo by the introduction of a vertically-shifted mesh, whose vertical levels are halfway between those of the primary mesh.
The top and bottom surfaces of the primary mesh and the vertically-shifted mesh coincide.
The density of a moisture species $\widetilde{\rho}_r$ is defined on this vertically-shifted mesh, using the same elements as $\mathbb{V}_\rho$ (with DoFs in cell centres), with this new space written as $\widetilde{\mathbb{V}}_\rho$, where the tilde $\widetilde{\cdot}$ denotes a quantity on the vertically-shifted mesh.
The vertically-shifted mesh then has one more level than the primary mesh, so that $\widetilde{\mathbb{V}}_\rho$ has the same number of DoFs as $\mathbb{V}_\theta$.
A similar mesh is used by \citet{thuburn2022numerical} to obtain entropy conservation with a Charney-Phillips staggering. \\
\\
The moisture density is calculated from $m_r$ and $\rho_d$ by converting the two fields to the $\widetilde{\mathbb{V}}_\rho$ space.
This uses two operators, $\mathcal{M}:\mathbb{V}_\theta\to\widetilde{\mathbb{V}}_\rho$ and $\mathcal{Q}:\mathbb{V}_\rho\to\widetilde{\mathbb{V}}_\rho$, so that
\begin{equation}
\widetilde{\rho}_r = \mathcal{M}[m_r]\times\mathcal{Q}[\rho_d],
\end{equation}
with the values of $\widetilde{\rho}_r$ given by the pointwise product of $\mathcal{M}[m_r]$ and $\mathcal{Q}[\rho_d]$.
The details of these operators will be discussed in Section \ref{sec:formulation_m}.
The dynamical core then conserves the following definition of moist mass:
\begin{equation} \label{eqn:moist mass}
\int_\varOmega \widetilde{\rho}_r \dx{V}.
\end{equation}
There is a vertically-shifted mesh corresponding to each mesh with different horizontal resolution, and so the shifting operators $\mathcal{M}$ and $\mathcal{Q}$ can also be defined on meshes with finer and coarser horizontal resolutions.
\begin{figure}[h!] 
\centering
\includegraphics[width=0.7\textwidth]{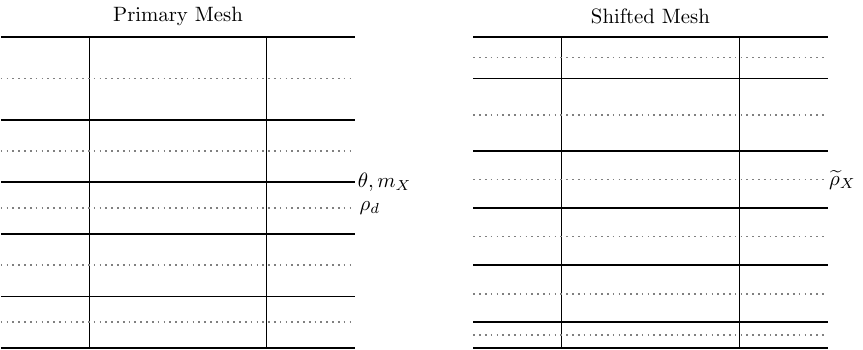}
\caption{A vertical cross-section illustrating the vertically-shifted mesh used in GungHo to describe moisture density, with solid black lines showing the top/bottom surfaces of elements and dotted grey lines showing the vertical centres of the levels.
The moisture mixing ratio $m_r$ is co-located with $\theta$ at the top and bottom surfaces of elements on the primary mesh, while the moisture density $\widetilde{\rho}_r$ is described at cell centres on a vertically-shifted mesh.
The top and bottom surfaces of elements on the vertically-shifted mesh coincide with the cell centres of elements on the primary mesh, so that the elements are shifted relative to those on the primary mesh.
The vertically-shifted mesh has one more level than the primary mesh.
} \label{fig:shifted}
\end{figure}

\subsection{Notation} \label{sec:notation}
\begin{figure}[h!] 
\centering
\includegraphics[width=0.5\textwidth]{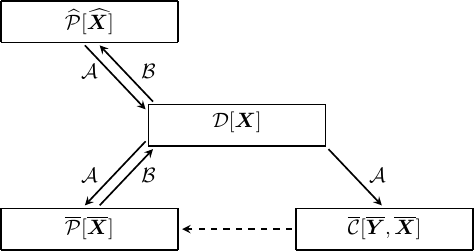}
\caption{A representation of a general atmospheric model with different components on different meshes, showing the three configurations considered in this work.
The dynamical core, described by operator $\mathcal{D}$, evolves the prognostic variables $\bm{X}$.
This is coupled to physical parametrisations $\widehat{\mathcal{P}}$ and $\overline{\mathcal{P}}$, which are computed on finer and coarser meshes respectively.
The final model component is $\overline{\mathcal{C}}$, which describes the evolution of $\overline{\bm{Y}}$, the chemical and aerosol variables.
These chemicals and aerosols may be used as auxiliary variables by a physics scheme (for instance a radiation scheme).
} \label{fig:model_schematic}
\end{figure}

It is convenient at this point to introduce the notation that is used in the rest of the paper.
Let the dynamical prognostic variables $\bm{X}$ evolved by the model be contained in some abstract space $\mathbb{V}_X$ so that $\bm{X}\in\mathbb{V}_X$, while the prognostic chemicals and aerosols $\bm{Y}$ are contained in a space $\mathbb{V}_Y$. \\
\\
These components may use meshes of different resolutions to one another.
Entities on a mesh that is finer resolution than that of the dynamical core are denoted with a hat $\widehat{\cdot}$ .
An overline $\overline{\cdot}$ denotes entities on a coarser mesh than that of the dynamical core.
As mentioned in the previous section, a tilde $\widetilde{\cdot}$ is used to denote entities on a vertically-shifted mesh.
Unadorned entities are on the same mesh as that used by the dynamical core.
\\
\\
With this notation, the components of the model described in Section \ref{sec:scope} that we will use in the remainder of the paper can be represented by the following operators:
\begin{enumerate}[leftmargin=*]
\item the dynamical core, $\mathcal{D}:\mathbb{V}_X\to\mathbb{V}_X$;
\item physics schemes that are computed on a finer mesh than the dynamical core, $\widehat{\mathcal{P}}:\widehat{\mathbb{V}}_X\to\widehat{\mathbb{V}}_X$;
\item physics schemes that are computed on a coarser mesh than the dynamical core, $\overline{\mathcal{P}}:\overline{\mathbb{V}}_X\to\overline{\mathbb{V}}_X$;
\item the chemistry and aerosol component on a coarser mesh than the dynamical core, $\overline{\mathcal{C}}:\left(\overline{\mathbb{V}}_Y,\overline{\mathbb{V}}_X\right)\to\overline{\mathbb{V}}_Y$.
\end{enumerate}
The interactions between components that use different meshes involve mapping fields from one mesh to another.
These mappings can also be represented by the action of operators:
\begin{equation}
\mathcal{A}:\widehat{\mathbb{V}}_X\to\mathbb{V}_X, \quad 
\mathcal{B}:\mathbb{V}_X\to\widehat{\mathbb{V}}_X,
\end{equation}
so that $\mathcal{A}$ maps fields to a coarser mesh, while $\mathcal{B}$ maps fields to a finer mesh.
Related operators can be defined to mapping fields between $\mathbb{V}_X$ and $\overline{\mathbb{V}}_X$, although for brevity these are also denoted by $\mathcal{A}$ and $\mathcal{B}$.
Thus $\mathcal{A}$ is akin to the \textit{restriction} operators used in the geometric multi-grid solver technique (see for instance \cite{maynard2020multigrid}, whereas $\mathcal{B}$ performs the role of a \textit{prolongation} operator.
It is also helpful to introduce the identification and reconstruction operators for mapping fields to finer meshes:
\begin{equation}
\mathcal{I}:\mathbb{V}_X\to\widehat{\mathbb{V}}_X, \quad \mathcal{R}:\mathbb{V}_X\to\widehat{\mathbb{V}}_X.
\end{equation}
In the absence of the orography (described in the next section) $\overline{\mathbb{V}}_X\subset\mathbb{V}_X\subset\widehat{\mathbb{V}}_X$. Fields on a coarser mesh can therefore be exactly represented, or \textit{identified}, on a finer mesh, with this operation denoted by $\mathcal{I}$.
The identification operators only use information from a single coarse cell to determine the value of a field in a cell on a finer mesh. 
In contrast, the \textit{reconstruction} operator $\mathcal{R}$ uses a stencil that takes field values from neighbouring coarse cells to obtain a higher-order reconstruction of the field.
These operators are discussed in more detail in Section \ref{sec:formulation}. \\  
\\
Some of the operators and their interactions are illustrated in Figure \ref{fig:model_schematic}, while are the operators are listed in Table \ref{tab:operators}.
\begin{table}[]
\centering
\begin{tabular}{l|c|c|c}
Operator & Notation & Domain & Co-domain \\
\hline
Dynamical core & $\mathcal{D}$ & $\mathbb{V}_X$ & $\mathbb{V}_X$ \\
Fine physics scheme &  $\widehat{\mathcal{P}}$ & $\widehat{\mathbb{V}}_X$ & $\widehat{\mathbb{V}}_X$ \\
Coarse physics scheme &  $\overline{\mathcal{P}}$ & $\overline{\mathbb{V}}_X$ & $\overline{\mathbb{V}}_X$ \\
Chemistry/aerosol model &  $\overline{\mathcal{C}}$ & $\left(\overline{\mathbb{V}}_Y,\overline{\mathbb{V}}_X\right)$ & $\overline{\mathbb{V}}_Y$ \\
Restriction & $\mathcal{A}$ & $\widehat{\mathbb{V}}_X$ & $\mathbb{V}_X$  \\
Prolongation & $\mathcal{B}$ & $\mathbb{V}_X$ & $\widehat{\mathbb{V}}_X$  \\
Identification & $\mathcal{I}$ & $\mathbb{V}_X$ & $\widehat{\mathbb{V}}_X$  \\
Reconstruction & $\mathcal{R}$ & $\mathbb{V}_X$ & $\widehat{\mathbb{V}}_X$  \\
Shifting operator for density & $\mathcal{Q}$ & $\mathbb{V}_\rho$ & $\widetilde{\mathbb{V}}_\rho$  \\
Shifting operator for mixing ratio & $\mathcal{M}$ & $\mathbb{V}_\theta$ & $\widetilde{\mathbb{V}}_\rho$
\end{tabular}
\caption{A list of the operators used in the formulation of Section \ref{sec:formulation}, showing the domain and co-domain.}
\label{tab:operators}
\end{table}

\subsection{Orography} \label{sec:orography}
\begin{figure}[h!] 
\centering
\includegraphics[width=0.5\textwidth]{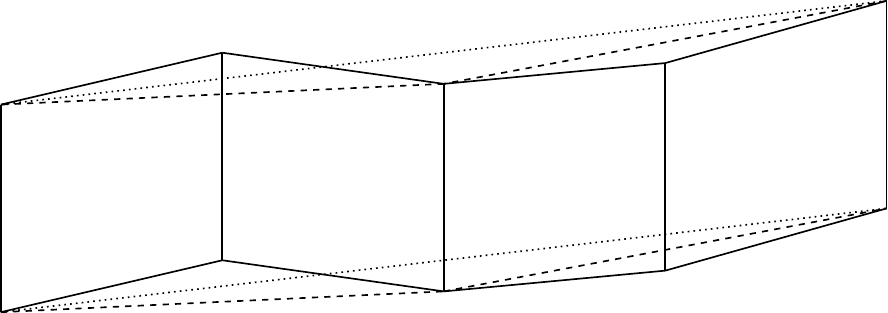}
\caption{An illustration of the strategy for describing the domain's orography for different meshes, through a vertical cross-section of one layer of elements.
The solid lines represent elements from the finest mesh, while dashed lines represent a mesh with intermediate resolution and the dotted lines showing the shape of the coarsest mesh.
The discretisation uses terrain-following coordinates, so the mesh's vertical coordinates are distorted to describe the orography.
The cell vertices of any coarser mesh are chosen to coincide with the appropriate vertices on the finest mesh, which defines the representation of the orography on the coarser meshes.} \label{fig:orography}
\end{figure}

GungHo uses terrain-following coordinates to describe the orography, so that the vertical coordinates of the mesh's vertices are modified to capture the effect of the planet's surface.
In general, the top and bottom faces of cells are sloped, while the lateral faces are aligned with the model's vertical direction.
When the model uses multiple meshes, the orography is first defined through the coordinates of the vertices on the finest mesh.
The vertices of cells in the coarser meshes are chosen to be coincident with the corresponding vertices on the finer mesh.
It should be noted that once the meshes have been modified to describe orography, cells in one layer on one mesh may overlap with cells of a different layer from another mesh.
The volume of the fine cells nested within a coarse cell may not necessarily equal the volume of the coarse cell.
This strategy is illustrated in Figure \ref{fig:orography}.
\section{Properties of Formulation} \label{sec:desirable properties}
Following the approach of \cite{herrington2019exploring} and \cite{herrington2019physics}, before introducing our formulation for coupling the components across different meshes, we list properties that we consider desirable for the formulation to possess.
Throughout Sections \ref{sec:desirable properties} and \ref{sec:formulation}, the properties will generally be discussed for mapping between $\mathbb{V}_X$ and $\widehat{\mathbb{V}}_X$, as the same operators are used for mapping between $\mathbb{V}_X$ and $\overline{\mathbb{V}}_X$.
\begin{enumerate}[leftmargin=*]
\item \textbf{Reversibility}. \label{prop:reversibility}
The combination of restriction and prolongation operators must be chosen so that mapping a field from a coarser mesh to a finer mesh and back results in an unchanged field, i.e.
\begin{equation}
\mathcal{A}\left[\mathcal{B}\left[\bm{X}\right]\right]=\bm{X}.
\end{equation}
This does not hold if the roles of $\mathcal{A}$ and $\mathcal{B}$ are reversed,
as information is lost as a field on a finer mesh is restricted to a coarser mesh.

\item \textbf{Preservation of a steady-state}. \label{prop:steady-state}
Consider a physical parametrisation that is computed upon a different mesh to the dynamical core.
If this physical parametrisation does not change the prognostic variables \textit{on the mesh of the physical parametrisation}, then the prognostic variables on the mesh of the dynamical core must not be changed by the combined process of mapping the prognostic fields to the physical parametrisation, computing the physical parametrisation and then mapping back.

\item \textbf{Conservation of mass of chemicals and aerosols}. \label{prop:chemical conservation}
When chemicals and aerosols are transported on the same mesh as the dynamical core, the masses of chemicals and aerosols are conserved.
This should still be true if these chemicals and aerosols are represented on a coarser mesh than the dynamical core, so that the transport of chemicals and aerosols conserves
\begin{equation}
\int_{\overline{\varOmega}} \mathcal{A}\left[\rho_d\right] \overline{a}_Y\dx{V},
\end{equation}
where $\overline{\varOmega}$ is the domain described by the coarser mesh.
\item \textbf{Preservation of constant chemical and aerosol mixing ratios}. \label{prop:chemical consistency}
The transport of chemicals and aerosols on a coarse mesh must preserve a constant mixing ratio.
This can be described as \textit{consistent transport}, as it implies that the chemical/aerosol densities evolve consistently with the density of dry air.

\item \textbf{Local conservation of mass of moisture species}. \label{prop:moisture conservation}
The dynamical core and physical parametrisations conserve the mass of moisture, in the absence of physical sources and sinks.
This conservation is local, in the sense that there is a local closed mass budget, as moisture obeys a conservative form of the transport equation.
The mapping operators for moisture should also locally conserve the mass of moisture locally within a coarse cell and over the fine cells contained within it.

\item \textbf{Preservation of constant mixing ratios of moisture species}. \label{prop:moisture consistency}
If a mixing ratio field takes a constant value $C$, then this must be preserved by the mixing ratio mapping operators (denoted by subscript $m$), so that
\begin{equation}
\mathcal{A}_m\left[C\right] = C, \quad \mathrm{and} \quad
\mathcal{B}_m\left[C\right] = C.
\end{equation}

\item \textbf{Avoid generation of negative moisture mixing ratios}. \label{prop:moisture positivity}
Negative values of moisture mixing ratios are unphysical and so must not be generated by the mapping formulation.
This is a weaker requirement than local shape preservation, which was considered by \citet{herrington2019exploring}, because the physical parametrisations themselves do not enforce local shape preservation, whereas they do ensure that negative values are not generated.

\item \textbf{Preservation of linear correlation of moisture mixing ratios}. \label{prop:moisture correlation}
If two moisture mixing ratios are linearly correlated on one mesh, so that $m_1=\alpha m_2 + \beta$ for constants $\alpha$ and $\beta$, then this linear correlation should hold after the two fields are mapped to another mesh.
As described by \citet{lauritzen2012evaluating}, these correlations can be important for determining the evolution of these variables.
This is also a property held by the approach of \cite{herrington2019exploring}.

\item \textbf{Accuracy}. \label{prop:accuracy}
The order of accuracy of the prolongation mapping should match the accuracy of the dynamical core.
For GungHo, this means second-order accuracy in space so that a field varying linearly in space should be exactly represented.
\end{enumerate}
As discussed by \citet{herrington2019physics}, conservation of other properties such as axial angular momentum, entropy or energy may be desirable but can be difficult to attain.
However, GungHo does not inherently conserve these properties so we do not see it as essential that they should be conserved by the formulation presented in the next section.

\section{Formulation} \label{sec:formulation}
To satisfy the desirable properties listed in Section \ref{sec:desirable properties}, we place two requirements on the operators in the formulation:
\begin{requirement} \label{req:simple reversibility}
The restriction operator $\mathcal{A}$ must act as the inverse of the identification operator $\mathcal{I}$, so that for any prognostic variable $\bm{X}$,
\begin{equation}
\mathcal{A}\left[\mathcal{I}\left[\bm{X}\right]\right] = \bm{X}.
\end{equation}
\end{requirement}
\begin{requirement}  \label{req:zero preservation}
The restriction operator $\mathcal{A}$ and the prolongation operator $\mathcal{B}$ must preserve a constant zero field, $\bm{0}$:
\begin{equation}
\mathcal{A}\left[\bm{0}\right]=\bm{0} \quad \mathrm{and} \quad \mathcal{B}\left[\bm{0}\right] = \bm{0}.
\end{equation}
\end{requirement}
\noindent Note that Requirement \ref{req:zero preservation} applies to all fields, while the stronger constraint of Property \ref{prop:moisture consistency} applies to just moisture mixing ratios.
\noindent Before discussing the restriction and prolongation operators for each of the prognostic variables, it is helpful to present features that are common to the operators for each of the scalar prognostic variables (the wind field is treated separately).
To obtain the reversibility discussed in Property \ref{prop:reversibility}, the prolongation operators are chosen for all scalar variables (with an additional subtlety for the moisture variables discussed in Section \ref{sec:formulation_m}) so that
\begin{equation} \label{def:B}
\mathcal{B}\left[\bm{X}\right] \equiv \mathcal{R}\left[\bm{X}\right] - \mathcal{I}\left[\mathcal{A}\left[\mathcal{R}\left[\bm{X}\right]\right]\right] + \mathcal{I}\left[\bm{X}\right].
\end{equation}
This has the same form as the recovery operator used by \citet{bendall2019recovered} and \citet{bendall2023improving} to obtain reversibility and mass conservation when recovering fields from lower to higher-order finite element spaces.
Then, given Requirement \ref{req:simple reversibility}, it can be seen that this structure for $\mathcal{B}$ will satisfy Property \ref{prop:reversibility}, as
\begin{equation} \label{eqn:general reversibility proof}
\mathcal{A}\left[\mathcal{B}\left[\bm{X}\right]\right]=
\mathcal{A}\left[\mathcal{R}\left[\bm{X}\right]\right]
-\mathcal{A}\left[\mathcal{I}\left[\mathcal{A}\left[\mathcal{R}\left[\bm{X}\right]\right]\right]\right]
+\mathcal{A}\left[\mathcal{I}\left[\bm{X}\right]\right]
= \mathcal{A}\left[\mathcal{R}\left[\bm{X}\right]\right]-\mathcal{A}\left[\mathcal{R}\left[\bm{X}\right]\right] + \bm{X}
= \bm{X},
\end{equation}
and the choice of \eqref{def:B} ensures that Property \ref{prop:reversibility} is obtained.
With the form of \eqref{def:B}, the reconstruction operator $\mathcal{R}$ defines the accuracy of the prolongation operator, while the remaining two terms can be considered as a correction to provide reversibility.
To meet Property \ref{prop:accuracy}, $\mathcal{R}$ should then be chosen to have the same order of accuracy as the dynamical core.
This form also means that the extrema of $\bm{X}$ will always lie within the extrema of $\mathcal{B}\left[\bm{X}\right]$.
\\
\\
To obtain Property \ref{prop:steady-state}, we take the same approach as \citet{herrington2019exploring}.
Denoting the field before and after the physical parametrisation by superscripts $n$ and $n+1$, so that $\bm{X}^{n+1}=\mathcal{P}\left[\bm{X}^n\right]$, then the increment corresponding the the physical parametrisation is simply
\begin{equation}
\Delta \mathcal{P}\left[\bm{X}^n\right] = \bm{X}^{n+1} - \bm{X}^n.
\end{equation}
To perform a physical parametrisation on a different mesh to the dynamical core, the updated prognostic fields are computed through
\begin{equation} \label{eqn:map_tendencies}
\bm{X}^{n+1} = \bm{X}^n + \mathcal{A}\left[\Delta\widehat{\mathcal{P}}\left[\mathcal{B}\left[\bm{X}^n\right]\right]\right],
\quad \mathrm{or} \quad
\bm{X}^{n+1} = \bm{X}^n + \mathcal{B}\left[\Delta\overline{\mathcal{P}}\left[\mathcal{A}\left[\bm{X}^n\right]\right]\right].
\end{equation}
Thus before physical parametrisations, prognostic variables are mapped from one mesh to another,
while after physical parametrisations, increments are mapped between meshes.
The situations considered by Property \ref{prop:steady-state} can be expressed in terms of increments, as if $\overline{\mathcal{P}}\left[\mathcal{A}\left[\bm{X}\right]\right] = \mathcal{A}\left[\bm{X}\right]$ then $\Delta\overline{\mathcal{P}}\left[\mathcal{A}\left[\bm{X}\right]\right]=\bm{0}$.
Provided that Requirement \ref{req:zero preservation} holds, then in this situation,
$\bm{X}^{n+1}=\bm{X}^n+\mathcal{B}\left[\bm{0}\right]=\bm{X}^n$.
A similar relation holds if the physical parametrisation is performed on a finer mesh.
\\
\\
However, the construction of \eqref{eqn:map_tendencies} makes satisfying the preservation of moisture positivity (Property \ref{prop:moisture positivity}) challenging when the physical parametrisation is computed on a coarser mesh. Although it is assumed that physical parametrisations do not generate negative moisture mixing ratio values on the mesh upon which they act, when the increment is mapped to the dynamical core mesh and added to the original mixing ratio field this can still generate spurious negative values.
The solution to this is discussed in Section \ref{sec:formulation_m}.
\\
\\
The remainder of this section specifies the particular restriction and prolongation operator for each prognostic variable, with a subscript to the operator denoting the variable, e.g.\ $\mathcal{A}_u$ for the restriction operator for the velocity $\bm{u}$.
The prolongation operators $\mathcal{B}_\rho$, $\mathcal{B}_\theta$ and $\mathcal{B}_\Pi$ take the form of \eqref{def:B}, so only the identification and restriction operators $\mathcal{I}$ and $\mathcal{A}$ need specifying.

\subsection{Mapping operators for the pressure and potential temperature fields} \label{sec:pressure mapping}
The mapping operators for the Exner pressure $\Pi$ and potential temperature $\theta$ are very similar.
The only difference is that $\Pi$ is expressed at points located in cell centres, while $\theta$ is vertically staggered from this.
As the vertical structure of the different meshes is the same, the operators involve only horizontal reconstruction or averaging.
Since the properties in Section \ref{sec:desirable properties} relating to $\Pi$ and $\theta$ are the same, the operators for $\Pi$ and $\theta$ take the same form as one another.
Therefore this section only presents the operators for $\Pi$. \\
\\
The restriction of $\Pi$ from a fine mesh to a coarse mesh consists of taking the arithmetic mean of the values in the fine cells contained within each coarse cell.
Let the Exner pressure field in the $j$-th fine cell within the $i$-th coarse cell in the $k$-th layer be denoted by $\widehat{\Pi}|_{i,j}^k$, and the value in the corresponding coarse cell be $\left.\Pi\right|_i^k$.
If there are $N_j$ fine cells in the $i$-th coarse cell then the action of $\mathcal{A}_\Pi$ is given by
\begin{equation}
\mathcal{A}_\Pi\left[\widehat{\Pi}\right] \equiv \left.\Pi\right|_i^k  = \frac{1}{N_j}\sum_{j=1}^{N_j} \widehat{\Pi}|_{i,j}^k.
\end{equation}
The identification operator $\mathcal{I}_\Pi$ is simply: 
\begin{equation}
\mathcal{I}_\Pi\left[\Pi\right] \equiv \widehat{\Pi}|_{i,j}^k =\left.\Pi\right|_i^k.
\end{equation}
Then this combination of $\mathcal{I}_\Pi$ and $\mathcal{A}_\Pi$ satisfies Requirement \ref{req:simple reversibility}, as
\begin{equation} \label{eqn:pressure reversibility}
\mathcal{A}_\Pi\left[\mathcal{I}_\Pi\left[\left.\Pi\right|_i^k\right]\right]=
\frac{1}{N_j}\sum_{j=1}^{N_j}\left.\Pi\right|_{i}^k = \left.\Pi\right|_i^k.
\end{equation}
The final operator is the reconstruction operator $\mathcal{R}_\Pi$, which uses a stencil over the $N_l$ neighbouring cells, with these cells indexed by $l$.
The operator is a simple weighted sum,
\begin{equation}
\mathcal{R}_\Pi\left[\Pi\right] \equiv \widehat{\Pi}|_{i,j}^k = \sum_{l=1}^{N_l} c_{i,j}^l \left.\Pi\right|_i^{k,l},
\end{equation}
where the coefficients $c^l_{i,j}$ sum to unity and can be chosen to give any particular reconstruction.
To give an order of accuracy approaching second-order, in this work the coefficients correspond to a linear reconstruction.
\\
\\
The operator $\mathcal{B}_\Pi$ can then be found from \eqref{def:B} to be described by
\begin{equation}
\mathcal{B}_\Pi\left[\Pi\right] \equiv \widehat{\Pi}|_{i,j}^k = \left.\Pi\right|_i^k + \sum_{l=1}^{N_l} c_{i,j}^l \left.\Pi\right|_i^{k,l} - \frac{1}{N_j}\sum_{m=1}^{N_j}\sum_{l=1}^{N_l} c_{i,m}^l \left.\Pi\right|_i^{k,l}.
\end{equation}
With these choices of operator, both $\mathcal{A}$ and $\mathcal{B}$ preserve a constant Exner pressure field or potential temperature fields (and hence also satisfy Requirement \ref{req:zero preservation}).

\subsection{Mapping operators for the density field}  \label{sec:rho mapping}
Key to achieving the local conservation of mass of moisture, chemical and aerosol species is choosing the dry density mapping operators so that they conserve mass within a coarse cell.
As discussed in Section \ref{sec:prognostics}, in GungHo density fields are represented by values at cell centres, and the basis functions for these fields are constant within a cell.
Then the mass in a cell is simply given by the product of the value of the density field for that cell with the cell's volume. \\
\\
Let the $k$-th cell in the $i$-th column be denoted by $e_i^k$, while the $j$-th cell on a finer mesh that is nested within it is $\widehat{e}_{i,j}^k$.
The restriction operator $\mathcal{A}_\rho\left[\widehat{\rho}\right]$ is defined by
\begin{equation}  \label{eqn:restriction_rho}
\mathcal{A}_\rho[\widehat{\rho}] \equiv \left.\rho\right|^k_i =
\frac{1}{\int_{e_i^k} \mathrm{d}V} \sum_{j=1}^{N_j} \left.\widehat{\rho}\right|^k_{i,j} \int_{\widehat{e}_{i,j}^k} \mathrm{d}V.
\end{equation}
This ensures that mass is conserved within a coarse cell by the restriction process.
If $\sum_{j=1}^{N_j}\int_{\widehat{e}^k_{i,j}}\mathrm{d}V=\int_{e^k_i}\mathrm{d}V$ then a constant density field is preserved by this restriction, but as discussed in Section \ref{sec:orography}, this is not necessarily true when the mesh is distorted by orography.
If the volume of the domain is different between the two meshes then it is not possible to both locally conserve mass and preserve a constant density. \\
\\
The identification operator $\mathcal{I}_\rho$ must conserve mass within a coarse element, so that it is given by
\begin{equation} \label{eqn:identification_rho}
\mathcal{I}_\rho\left[\rho\right]\equiv \left.\widehat{\rho}\right|_{i,j}^k =
\frac{\int_{e_i^k} \mathrm{d}V}{N_j \int_{\widehat{e}_{i,j}^k}\mathrm{d}V} \left.\rho\right|_i^k , 
\end{equation}
which also combines with $\mathcal{A}_\rho$ to satisfy Requirement \ref{req:simple reversibility}. \\
\\
The reconstruction operator $\mathcal{R}_\rho$ does not need to conserve mass, as conservation of mass is only required by $\mathcal{B}_\rho$.
Therefore the reconstruction operator $\mathcal{R}_\rho$ is taken to be $\mathcal{R}_\Pi$.
Conservation of mass of $\mathcal{B}_\rho$ follows from conservation of mass of $\mathcal{A}_\rho$ and $\mathcal{I}_\rho$.

\subsection{Mapping operators for the wind field}  \label{sec:wind mapping}
In the formulation considered in this work, the physical parametrisations that provide increments to the wind are not computed on a different mesh to the dynamical core mesh, although scalar quantities contributing to those physical parametrisations may be calculated on different meshes.
However the wind field must still be mapped to other meshes, for the transport of aerosols on a coarser mesh and also as an auxiliary field to physical parametrisations that return increments to other scalar fields. \\
\\
The wind field is described in GungHo through its normal component to cell faces.
We describe the faces of fine cells that coincide with the faces of coarse cells as being on the \textit{exterior} of coarse cells, while those faces that do not coincide are on the \textit{interior} of coarse cells.
The operators presented in this section are motivated by the consistent transport of chemicals and aerosols which will be discussed in Section \ref{sec:consistent transport}.
Key to this is for the operators to conserve the velocity flux through the faces of coarse cells.\\
\\
Let the $N_f$ faces of the element $e_i^k$ be denoted by $\Gamma_{i,f}^k$ (so that faces are indexed by $f$).
Similarly, the faces of the element $\widehat{e}_{i,j}^k$ are given by $\widehat{\Gamma}_{i,j,f}^k$.
However the face $\Gamma_{i,f}^k$ coincides with $N_g$ faces of fine elements, which can be written as $\widehat{\Gamma}_{i,f}^{k,g}$, where $g$ is the index of the coincident fine faces.
The value of $N_g$ may be different for different faces $\Gamma_{i,f}^k$.
The variables $\left.u\right|_{i,f}^k$ and $\left.\widehat{u}\right|_{i,f}^{k,g}$ are the contravariant wind components that correspond to the faces $\Gamma_{i,f}^k$ and $\widehat{\Gamma}_{i,j,f}^k$. \\
\\
With this notation, the restriction operator $\mathcal{A}_u$ is defined through
\begin{equation} \label{eqn:restriction_u}
\mathcal{A}_u\left[\widehat{\bm{u}}\right] \equiv \left.u\right|_{i,f}^k = 
\frac{1}{\int_{\Gamma_{i,f}^k} \mathrm{d}A} \sum_{g=1}^{N_g} \left.\widehat{u}\right|_{i,f}^{k,g} \int_{\widehat{\Gamma}_{i,f}^{k,g}} \mathrm{d}A,
 \end{equation}
where $\mathrm{d}A$ is the measure of the surface integral for a cell face.
Only those fine mesh values that are on the exterior of coarse cells contribute to the restriction. \\
\\
The prolongation operator $\mathcal{B}_u$ takes a different form to those used for the scalar fields, as identification and reconstruction operators do not get defined.
The fine cell values on the exterior of coarse cells are obtained through 
\begin{equation}
\mathcal{B}_u[\bm{u}]=\left.\widehat{u}\right|^k_{i,f} = 
\frac{\int_{\Gamma_{i,f}^k} \mathrm{d}A}{N_g \int_{\widehat{\Gamma}_{i,f}^{k,g}} \mathrm{d}A} \left.u\right|_{i,f}^k.
\end{equation}
The horizontal wind for the faces of fine cells that are interior to coarse cells are obtained through linear interpolation of the values from opposite faces of the coarse cell.
As the $\mathbb{V}_u$ basis functions are linear functions in the direction of the normal component, this prolongation emulates an identification operator.
\\
\\
With these choices of operator, Requirement \ref{req:zero preservation} is satisfied as the zero vector is mapped from one mesh to another.
Since the wind values for the faces on the interior of the coarse cell do not contribute to the restriction operator, these do not need to be considered.
As the wind field does not directly have physics increments computed on different meshes, it is not necessary to build a higher-order reconstruction operator $\mathcal{R}_u$.
\subsection{Conservative and consistent transport of chemicals and aerosols} \label{sec:consistent transport}
Local mass conservation of tracers is achieved by transporting the density $\rho_Y$ using a conservative form of the transport equation.
If the mass fluxes of dry air and of the $Y$-th tracer species are defined as $\bm{F}_d:=\rho_d\bm{v}$ and $\bm{F}_Y:=\rho_Y\bm{v}$, equations \eqref{eqn:compressible_euler_rho} and \eqref{eqn:compressible_euler_tracer} for the transport of dry density and tracers can be written as
\begin{equation}
\pfrac{\rho_d}{t} + \bm{\nabla}\bm{\cdot}\bm{F}_d = 0, \quad
\pfrac{\rho_Y}{t} + \bm{\nabla}\bm{\cdot}\bm{F}_Y = 0,
\end{equation}
where the sources and sinks of the tracers have been omitted. 
Following the approach taken by \citet{lauritzen2011numerical,lauritzen2014standard,zangl2015icon} and \citet{thuburn2022numerical}, the tracer mass flux can be expressed as $\bm{F}_Y=m_Y\bm{F}_d$ such that the tracer transport obeys
\begin{equation} \label{eqn:consistent transport}
\pfrac{\rho_Y}{t}+\bm{\nabla\cdot}\left(m_Y\bm{F}_d\right) = 0.
\end{equation}
Using the same dry mass flux $\bm{F}_d$ to transport both $\rho_d$ and $\rho_Y$ is key to ensuring consistent tracer transport.
However as $\rho_Y$ is transported on a coarser mesh than $\rho_d$, it is necessary to map $\bm{F}_d$ to the coarser mesh.
The approach described in this section is similar to the framework presented by \citet{bendall2023solution} used for conservative and consistent transport of moisture species on a vertically-shifted mesh.\\
\\
To begin, it is assumed that the discretised transport of $\rho_d\in\mathbb{V}_\rho$ and $\overline{\rho}_Y\in\overline{\mathbb{V}}_\rho$ can be expressed as two-time level schemes:
\begin{equation} \label{eqn:tracers two-step}
\rho_d^{n+1}=\rho_d^n - \Delta t \bm{\nabla}\bm{\cdot}\bm{F}_d
\quad \mathrm{and} \quad
\overline{\rho}_Y^{n+1} = \overline{\rho}_Y^n - \Delta t \overline{\bm{\nabla}}\bm{\cdot}\mathcal{F}\left[\overline{a}_Y,\overline{\bm{F}}_d\right],
\end{equation}
with the superscript $n$ denoting a field at the $n$-th time level and where the flux $\overline{\bm{F}}_Y$ has been calculated by a flux operator $\mathcal{F}:\overline{\mathbb{V}}_\rho,\overline{\mathbb{V}}_u\to\overline{\mathbb{V}}_u$.
The divergence operators act such that $\bm{\nabla\cdot}:\mathbb{V}_u\to\mathbb{V}_\rho$ and $\overline{\bm{\nabla}}\bm{\cdot}:\overline{\mathbb{V}}_u\to\overline{\mathbb{V}}_\rho$.
These satisfy the divergence theorem within a cell, so that
\begin{equation} 
\int_{e_i^k}\bm{\nabla\cdot u}\dx{V}=\sum_{f=1}^{N_f}\int_{\Gamma_{i,f}^k} \bm{u}\bm{\cdot}\mathrm{d}\bm{A},
\end{equation}
which becomes
\begin{equation} \label{eqn:div}
\left.(\bm{\nabla\cdot u})\right|_i^k \int_{e_i^k}\mathrm{d}V = \sum_{f=1}^{N_f} \left. u\right|_{i,f}^k\int_{\Gamma_{i,f}^k} \mathrm{d}A
\end{equation}
for the lowest-order finite elements used in GungHo.
As $\overline{a}_Y\in\overline{\mathbb{V}}_\rho$, the conversion between mixing ratio and density can be computed pointwise through
\begin{equation} \label{eqns:tracers conversions}
\overline{\rho}_Y = \overline{a}_Y\times \mathcal{A}_\rho\left[\rho_d\right]
\quad \mathrm{and} \quad
\overline{a}_Y = \overline{\rho}_Y / \mathcal{A}_\rho\left[\rho_d\right].
\end{equation}
Since $\overline{a}_Y$ is constant within a (coarse) cell and $\mathcal{A}_\rho\left[\rho_d\right]$ conserves mass within a coarse cell, the mass that is conserved by solving \eqref{eqn:consistent transport} is
\begin{equation}
\int_{\overline{\varOmega}} \overline{\rho}_Y\dx{V} = \int_{\overline{\varOmega}} \overline{a}_Y\mathcal{A}_\rho\left[\rho_d\right]\dx{V}.
\end{equation}
There are two further requirement concerning the operators in this set up:
\begin{requirement} \label{req:flux op}
The flux operator $\mathcal{F}$ must satisfy, for any constant $C\in\overline{\mathbb{V}}_\rho$ and any $\overline{\bm{F}}_d\in\overline{\mathbb{V}}_u$,
\begin{equation}
\mathcal{F}\left[C,\overline{\bm{F}}_d\right]=C\overline{\bm{F}}_d.
\end{equation}
\end{requirement}
\noindent As the focus here is on the operators to map between meshes, we do not detail a flux operator that meets Requirement \ref{req:flux op}, but one using a Method of Lines scheme with Runge-Kutta time stepping that does meet this requirement was presented by \citet{bendall2023solution}.
\begin{requirement} \label{req:consistent transport}
For all $\bm{F}_d\in\mathbb{V}_u$, the restriction and divergence operators commute so that
\begin{equation}
\mathcal{A}_\rho\left[\bm{\nabla}\bm{\cdot}\bm{F}_d\right] = \overline{\bm{\nabla}}\bm{\cdot}\mathcal{A}_u\left[\bm{F}_d\right].
\end{equation}
\end{requirement}
\noindent The combination of restriction operators \eqref{eqn:restriction_rho} and \eqref{eqn:restriction_u}  satisfy Requirement \ref{req:consistent transport} given the divergence operator \eqref{eqn:div}, as for $\bm{u}\in\mathbb{V}_u$, the operation $\mathcal{A}_\rho\left[\bm{\nabla\cdot u}\right]$ can be expressed through
\begin{subequations}
\begin{align}
\sum_{j=1}^{N_j} \left.(\bm{\nabla\cdot u})\right|_{i,j}^k \int_{e_{i,j}^k} \mathrm{d}V & =
\sum_{j=1}^{N_j}\sum_{f=1}^{N_f} \left.u\right|_{i,j,f}^k \int_{\Gamma_{i,j,f}^k}\mathrm{d}A,
\intertext{which is a sum over all the faces of the fine cells within the coarse cell.
However, for faces on the interior of a coarse cell, an integral over each face is exactly cancelled by another opposite integral, and so the final sum only includes those fluxes over the exterior faces of the coarse cell:}
& = \sum_{f=1}^{N_f} \sum_{g=1}^{N_g} \left.u\right|_{i,f}^{k,g} \int_{\Gamma_{i,f}^{k,g}} \mathrm{d}A,  \\
& = \sum_{f=1}^{N_f} \left.\overline{u}\right|_{i,f}^k \int_{\overline{\Gamma}_{i,f}^k} \mathrm{d}A,
\end{align}
\end{subequations}
which is equivalent to $\overline{\bm{\nabla}}\bm{\cdot}\mathcal{A}_u\left[\bm{u}\right]$. \\
\\
With these requirements, the consistent transport described in Property \ref{prop:chemical consistency} can be attained by a discretisation that uses \eqref{eqn:tracers two-step} and \eqref{eqns:tracers conversions}, as a constant mixing ratio field can be preserved.
Combining \eqref{eqn:tracers two-step} and \eqref{eqns:tracers conversions},
\begin{subequations}
\begin{align}
\overline{a}_Y^{n+1}
& = \overline{\rho}^{n+1}_Y / \mathcal{A}_\rho\left[\rho^{n+1}_d\right], \\
& = \left(\overline{\rho}_Y^n - \Delta t \overline{\bm{\nabla}}\bm{\cdot}\mathcal{F}\left[\overline{a}_Y,\overline{\bm{F}}_d\right] \right)
/ \mathcal{A}_\rho\left[\rho_d^n - \Delta t \bm{\nabla}\bm{\cdot}\bm{F}_d\right], \\
& = \left(\overline{a}_Y^n\mathcal{A}\left[\rho_d^n\right] - \Delta t \overline{\bm{\nabla}}\bm{\cdot}\mathcal{F}\left[\overline{a}_Y,\overline{\bm{F}}_d\right] \right)
/ \mathcal{A}_\rho\left[\rho_d^n - \Delta t \bm{\nabla}\bm{\cdot}\bm{F}_d\right].
\intertext{Inserting $\overline{a}_Y^n=C$ for a constant $C$ into the right-hand side and then using Requirements \ref{req:flux op} and \ref{req:consistent transport}, and that the restriction operator $\mathcal{A}_\rho$ is linear,}
& = \left(C\mathcal{A}\left[\rho_d^n\right] - \Delta t \overline{\bm{\nabla}}\bm{\cdot}\mathcal{F}\left[C,\overline{\bm{F}}_d\right] \right)
/ \mathcal{A}_\rho\left[\rho_d^n - \Delta t \bm{\nabla}\bm{\cdot}\bm{F}_d\right], \\
& = C
\end{align}
\end{subequations}
and a constant mixing ratio is preserved. \\
\\
Although in this section we have considered tracers in $\overline{\mathbb{V}}_\rho$, it is straightforward to extend this approach to transporting tracers in $\overline{\mathbb{V}}_\theta$ by using a vertically-shifted coarse mesh and restricting $\widetilde{\rho}_d$ and $\widetilde{\bm{F}}_d$ to this mesh.

\subsection{Mapping operators for the moisture mixing ratios}  \label{sec:formulation_m}
Before presenting the restriction and prolongation operators for the moisture mixing ratios, it is convenient to discuss the operators for converting $m_r\in\mathbb{V}_\theta$ to $\widetilde{\rho}_r\in\widetilde{\mathbb{V}}_\rho$ on the shifted mesh.
This involves the operators $\mathcal{M}$ and $\mathcal{Q}$, as described in Section \ref{sec:moisture_conservation}.
The forms of these operators are taken from \citet{bendall2023solution}, so that $\mathcal{Q}$ is defined by
\begin{equation} \label{eqn:shifted_mesh_operator}
\mathcal{Q}\left[\rho\right] \equiv \left.\widetilde{\rho}\right|^k_{i,j} = \frac{1}{2\int_{\widetilde{e}_{i,j}^k} \mathrm{d}V}\left[
\left.\rho\right|_{i,j}^k\int_{e_{i,j}^k}\mathrm{d}V+\left.\rho\right|_{i,j}^{k+1}\int_{e_{i,j}^{k+1}} \mathrm{d}V\right]
\quad \mathrm{for} \ k\in[2,N_k],
\end{equation}
if there are $N_k$ DoFs per column in $\mathbb{V}_\rho$, while in the top and bottom layers:
\begin{equation}
\left.\widetilde{\rho}\right|^1_{i,j} = \frac{\int_{e_{i,j}^{1}}\mathrm{d}V}{2\int_{\widetilde{e}_{i,j}^1} \mathrm{d}V}
\left.\rho\right|_{i,j}^1, \quad
\left.\widetilde{\rho}\right|^{N_k+1}_{i,j} = \frac{\int_{e_{i,j}^{N_k}}\mathrm{d}V}{2\int_{\widetilde{e}_{i,j}^{N_k+1}} \mathrm{d}V}
\left.\rho\right|_{i,j}^{N_k}.
\end{equation}
The operator $\mathcal{M}$ is defined by the pointwise assignment of values, with some interpolation in the top and bottom layers:
\begin{equation}
\left.\widetilde{m}_r\right|_i^1 =
\frac{1}{2}\left(\left.m_r\right|_i^1+\left.m_r\right|_i^2 \right), \quad
\left.\widetilde{m}_r\right|_i^{N_k+1} =
\frac{1}{2}\left(\left.m_r\right|_i^{N_k}+\left.m_r\right|_i^{N_k+1} \right), \quad
\left.\widetilde{m}_r\right|_i^{k}
= \left.m_r\right|_i^{k}
\quad \mathrm{for} \ k\in[2,N_k].
\end{equation}
The operator $\mathcal{M}^{-1}:\widetilde{\mathbb{V}}_\rho\to\mathbb{V}_\theta$ is the inverse of $\mathcal{M}$, and is straightforward to compute for this choice of $\mathcal{M}$.
Instead of interpolating, the values at the top and bottom of the column are obtained through linear extrapolation (with a correction to avoid the generation of negative values).
The operators $\mathcal{M}$ and $\mathcal{M}^{-1}$ preserve a constant mixing ratio field, and are both linear. \\
\\
\noindent Since moisture conservation is defined through the density $\widetilde{\rho}_r$, it is necessary to specify how the restriction and identification operators interact with the shifted mesh.
Using the definitions \eqref{eqn:restriction_rho}, \eqref{eqn:identification_rho} and \eqref{eqn:shifted_mesh_operator} of the operators $\mathcal{A}_\rho$, $\mathcal{I}_\rho$ and $\mathcal{Q}$, it can be shown that $\mathcal{Q}$ commutes with both $\mathcal{A}_\rho$ and $\mathcal{I}_\rho$ so that for any $\widehat{\rho}\in\widehat{\mathbb{V}}_\rho$ or $\rho\in\mathbb{V}_\rho$
\begin{equation}
\mathcal{Q}\left[\mathcal{A}_\rho\left[\widehat{\rho} \, \right]\right] =
\mathcal{A}_\rho\left[\mathcal{Q}\left[\widehat{\rho} \, \right]\right] 
\quad \mathrm{and} \quad
\mathcal{Q}\left[\mathcal{I}_\rho\left[\rho\right]\right] = 
\mathcal{I}_\rho\left[\mathcal{Q}\left[\rho\right]\right].
\end{equation}
\subsubsection{Restriction and Identification operators}
With these definitions, the restriction operator for the mixing ratio field is $\mathcal{A}_m$, which can be written in terms of existing operators as
\begin{equation}
\mathcal{A}_m\left[\widehat{m}_r\right]\equiv\mathcal{M}^{-1}\left[\mathcal{A}_\rho\left[\mathcal{M}\left[\widehat{m}_r\right]\times\mathcal{Q}\left[\widehat{\rho}_d\right]\right] / \mathcal{Q}\left[\mathcal{A_\rho}\left[\widehat{\rho}_d\right]\right]\right],
\end{equation}
while the related identification operator $\mathcal{I}$ is given by
\begin{equation}
\mathcal{I}_m\left[m_r\right] \equiv \mathcal{M}^{-1}\left[\mathcal{I}_\rho\left[\mathcal{M}\left[m_r\right]\times\mathcal{Q}\left[\rho_d\right] \right] / \mathcal{Q}\left[\mathcal{B}_\rho\left[\rho_d\right]\right]\right].
\end{equation}
These choices are designed so that Requirement \ref{req:simple reversibility} is satisfied: the operators involve expressing the moisture field as a density and then restricting or identifying that density.
As $\mathcal{A}_\rho\left[\mathcal{I}_\rho\left[\rho\right]\right] = \rho$, it then follows that $\mathcal{A}_m\left[\mathcal{I}_m\left[m_r\right]\right] = m_r$.\\
\\
By construction, these operators also provide Property \ref{prop:moisture conservation}, since the restriction and identification processes act upon a density field, so mass is naturally conserved by the mappings.
As all of the constituent operators are linear, then $\mathcal{A}_m$ and $\mathcal{I}_m$ are also linear and so satisfy Property \ref{prop:moisture correlation}.
The restriction operator $\mathcal{A}_m$ preserves a constant mixing ratio (Property \ref{prop:moisture consistency}) since
\begin{equation}
\mathcal{A}_m[C]
=\mathcal{M}^{-1}\left[\mathcal{A}_\rho\left[C\mathcal{Q}\left[\widehat{\rho}_d\right]\right] / \mathcal{Q}\left[\mathcal{A}_\rho\left[\widehat{\rho_d}\right]\right]\right]
=\mathcal{M}^{-1}\left[C\mathcal{A}_\rho\left[\mathcal{Q}\left[\widehat{\rho}_d\right]\right] / \mathcal{A}_\rho\left[\mathcal{Q}\left[\widehat{\rho_d}\right]\right]\right]=C,
\end{equation}
as $\mathcal{Q}$ commutes with $\mathcal{A}_\rho$ and $\mathcal{I}_\rho$. \\
\\
Finally, provided that $\mathcal{B}_\rho[\rho_d]$ is positive (which it should be for well-behaved density fields), then $\mathcal{I}_m$ and $\mathcal{A}_m$ cannot generate negative mixing ratio values,
as none of the operators $\mathcal{M}^{-1}$, $\mathcal{M}$, $\mathcal{Q}$, $\mathcal{A}_\rho$ and $\mathcal{I}_\rho$ can generate negative values. 
\subsubsection{Prolongation operator}
An initial prolongation operator is defined by
\begin{equation}
\mathcal{B}^\dagger_m\left[m_r\right] \equiv
\mathcal{R}_\theta\left[m_r\right] - \mathcal{I}_m\left[\mathcal{A}_m\left[\mathcal{R}_\theta\left[m_r\right]\right]\right]
+ \mathcal{I}_m\left[m_r\right].
\end{equation}
As this has the same structure as \eqref{def:B}, it satisfies Property \ref{prop:reversibility}.
Then a modified prolongation operator is
\begin{equation} \label{def:Bm}
\mathcal{B}_m\left[m_r\right] \equiv
(1-\lambda)\mathcal{R}_\theta\left[m_r\right] - (1-\lambda)\mathcal{I}_m\left[\mathcal{A}_m\left[\mathcal{R}_\theta\left[m_r\right]\right]\right]
+ \mathcal{I}_m\left[m_r\right]
\end{equation}
where $\lambda$ is a field in the same space on the same mesh as $m_r$, and takes values between 0 and 1.
The inclusion of $\lambda$ is to prevent the generation of negative mixing ratios, which will be discussed in Section \ref{sec:negativity}.
This operator can also be expressed as:
\begin{equation}
\mathcal{B}_m=(1-\lambda)\mathcal{B}_m^\dagger + \lambda\mathcal{I}_m.
\end{equation}
Since $\mathcal{A}_m$ is a linear operator, and $\mathcal{B}_m$ is a linear combination of $\mathcal{B}^\dagger_m$ and $\mathcal{I}_m$, $\mathcal{B}_m$ also satisfies Property \ref{prop:reversibility}.
As with $\mathcal{A}_m$ and $\mathcal{I}_m$, the operator $\mathcal{B}_m$ is linear and conserves mass within a coarse element.
A constant mixing ratio is also preserved by $\mathcal{B}^\dagger_m$ as
\begin{equation}
\mathcal{B}^\dagger_m\left[C\right] = \mathcal{R}_\theta\left[C\right] - \mathcal{I}_m\left[\mathcal{A}_m\left[\mathcal{R}_\theta\left[C\right]\right]\right] + \mathcal{I}_m\left[C\right] = C - \mathcal{I}_m\left[\mathcal{A}_m\left[C\right]\right] + \mathcal{I}_m\left[C\right] = C - \mathcal{I}_m\left[ C\right] + \mathcal{I}_m\left[C\right] = C,
\end{equation}
so $\mathcal{B}_m^\dagger$ preserves a constant.
The same steps can also be used to show that $\mathcal{B}_m$ preserves a constant mixing ratio.
\subsubsection{Prevention of negative mixing ratios} \label{sec:negativity}
\begin{figure}[h!]
\centering
\includegraphics[width=0.8\textwidth]{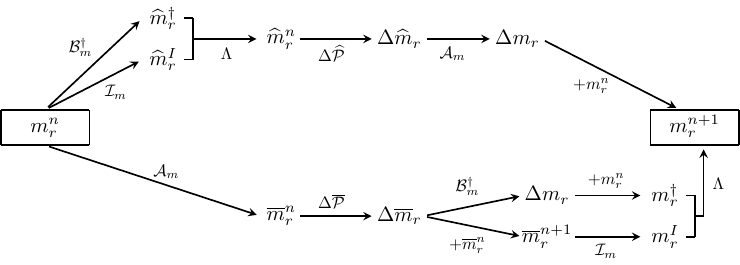}
\caption{The procedure to compute physical parametrisations for moisture mixing ratios upon finer or coarser meshes, including the steps to prevent the generation of negative values.
The mixing ratio field on the dynamical core mesh before the parametrisation is $m_r^n$, while the resulting mixing ratio field is given by $m_r^{n+1}$.
The upper half of the diagram describes a physical parametrisation on a finer mesh, where negative values can be generated by the prolongation to the finer mesh.
The lower half represents a physical parametrisation on a coarser mesh, where negative values could be caused by the addition of a tendency $\Delta \overline{m}_r$ computed on a coarse mesh to the original field.}
\label{fig:moisture_mapping}
\end{figure}

In this formulation, there are two situations in which negative moisture mixing ratios can be generated by the mapping process, unless care is taken.
The first is the prolongation of a mixing ratio field to a finer mesh.
The second is the addition of an increment to a mixing ratio field, when that increment has been calculated on a coarser mesh.
The solution in both situations involves combining the field which may be negative with one that is guaranteed not to be.
This is described through the operator $\Lambda:\mathbb{V}_\theta,\mathbb{V}_\theta\to\mathbb{V}_\theta$. \\
\\
To determine the action of $\Lambda$ and the value of $\lambda$, 
consider $m_r^{-}\in\mathbb{V}_\theta$, a mixing ratio field which may contain negative values, and
$m_r^{+}\in\mathbb{V}_\theta$, whose values are guaranteed not to be negative.
It is possible to define an operator $\Lambda$ which blends $m_r^{-}$ and $m_r^{+}$ to create a field $m_r$ which is also guaranteed not to be negative, through
\begin{equation} \label{eqn:Lambda}
m_r = \Lambda\left[m_r^{-},m_r^{+}\right] \equiv (1-\lambda)m_r^{-} + \lambda m_r^{+},
\end{equation}
where $\lambda\in\overline{\mathbb{V}}_\theta$ is a field on a coarser mesh whose values lie between 0 and 1. \\
\\
To find appropriate values of $\lambda$, consider the values of the mixing ratio field in one coarse cell.
If $m_r^-|_{i,j}^k<0$, then $m_r|_{i,j}^k=0$ if
\begin{equation}
\lambda|_i^k = \frac{-m_r^-|_{i,j}^k}{m_r^+|_{i,j}^k-m_r^-|_{i,j}^k},
\end{equation}
which is found by rearranging \eqref{eqn:Lambda}.
Negativity will be prevented by taking
\begin{equation} \label{eqn:positivity factor}
\left.\lambda\right|^k_i = \left\lbrace
\begin{matrix}
0, & \min_j\left(\left.m^{-}_r\right|_{i,j}^k\right) \geq 0 \\
\max_l\left(\dfrac{-\left.m^{-}_r\right|_{i,l}^k}
{\left.m_r^{+}\right|_{i,l}^k-\left.m_r^{-}\right|_{i,l}^k}\right), & \mathrm{otherwise}, \ \mathrm{where} \ l\in\curlybrackets{j:\left.m_r^-\right|_{i,j}^k<0}.
\end{matrix} \right.
\end{equation}
Since $\mathcal{I}_m$ does not generate negative mixing ratio values, the inclusion of $\lambda$ in the definition of $\mathcal{B}_m$ in \eqref{def:Bm} also prevents the generation of negative values, as
\begin{equation}
\mathcal{B}_m\left[m_r\right] = \Lambda\left[\mathcal{B}_m^\dagger\left[m_r\right], \mathcal{I}_m\left[m_r\right]\right] \equiv (1-\lambda)\mathcal{B}_m^\dagger\left[m_r\right] + \lambda \mathcal{I}_m\left[m_r\right].
\end{equation}
The positivity factor then acts like the slope limiter of \cite{barth1989design}, with the higher-order reconstruction obtained from $\mathcal{B}^\dagger_m$ limited by the minimum amount to ensure positivity.
A similar approach is used by \cite{herrington2019exploring} to tackle the same problem.
To ensure that linear correlations between moisture species are preserved, the same $\lambda$ field should be used for all species, which can be computed in each cell to be the maximum $\lambda$ value for each individual species. \\
\\
The other situation in which negative moisture values can be generated follows the computation of a physical parametrisation on a coarse mesh.
In this case, consider the mixing ratio $m^\dagger_r$ which has been updated with an increment from a coarse physical parametrisation, so that
\begin{equation}
m_r^\dagger = m^n_r + \mathcal{B}_m^\dagger\left[\Delta\overline{\mathcal{P}}\left[\mathcal{A}_m\left[m^n_r\right]\right]\right].
\end{equation}
This is not guaranteed to be positive.
However, assuming that the physical parametrisation is positivity-preserving and again using that $\mathcal{I}_m$ does not generate negative values, the following field is assured to be positive:
\begin{equation}
m_r^I:=\mathcal{I}_m\left[\mathcal{A}_m\left[m_r^n\right]+\Delta\overline{\mathcal{P}}\left[\mathcal{A}_m\left[m_r^n\right]\right]\right],
\end{equation}
which is the field (and not the increment) resulting from the physical parametrisation on the coarse mesh, mapped back to the original mesh using $\mathcal{I}_m$.
The resulting field $m_r^{n+1}$ on the dynamical core mesh is then computed from
\begin{equation}
m_r^{n+1} = \Lambda\left[m_r^\dagger, m_r^I\right],
\end{equation}
which avoids the generation of negative values.
The moisture mapping processes are summarised in Figure  \ref{fig:moisture_mapping}.

\section{Idealised Test Results} \label{sec:results}
The following section aims to demonstrate the formulation described in Section \ref{sec:formulation} through a series of idealised test cases.
While some tentative conclusions about the accuracy of these choices are highlighted, the primary motivation is to test the formulation without the complexities of a full suite of physical parametrisations.
The test cases are run using the GungHo dynamical core described in Section \ref{sec:lfric}, which as described in \citet{melvin2019mixed} and \citet{kent2023mixed}, uses an iterative semi-implicit time stepping scheme with a nested outer-inner loop structure like that of ENDGame, \cite{wood2014inherently}. Transport terms are treated explicitly in the outer loop using a Method of Lines (MoL) structure with finite volume spatial discretisation. Faster terms describing wave motions are treated implicitly in the inner loop, which consists of an iterative Newton solve. 
\\
\\
All variables are transported using the MoL scheme described in \citet{melvin2019mixed}, with vertical-horizontal Strang split to reduce the number of required substeps when the vertical Courant number is large. Dry density is transported conservatively, while the potential temperature and wind are transported in advective forms. When included, moisture species are transported conservatively and consistently with the scheme described in \citet{bendall2023solution}.
\\
\\
The cloud microphysics scheme, used in two of the test cases, is a simple evaporation-condensation scheme with latent heat feedback, so that the moisture species are water vapour and cloud liquid. The scheme is called after the transport step within the outer loop of the algorithm.

\subsection{Tracer transport on the sphere} \label{sec:spherical transport test}

\begin{figure}[h!]
\centering
\includegraphics[width=0.9\textwidth]{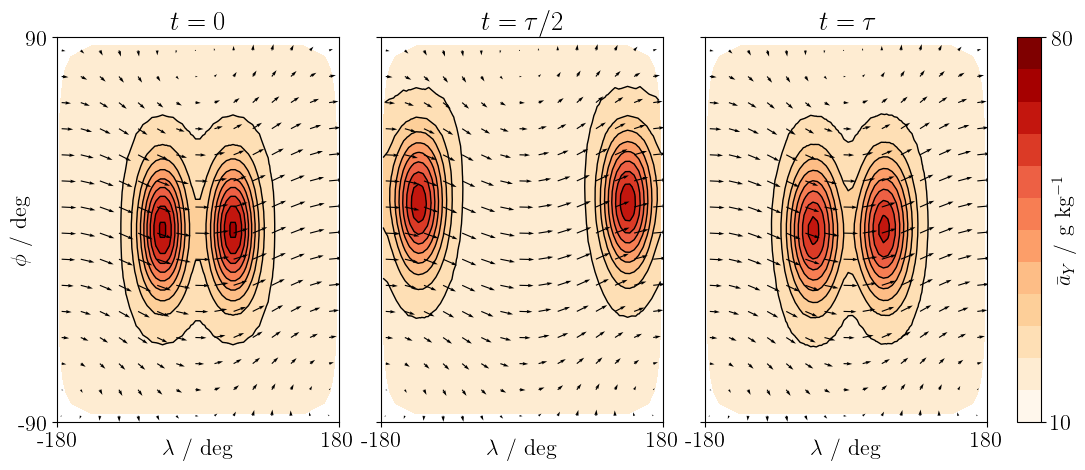}
\caption{The $\overline{a}_Y$ field used in the transport test case on the surface of a sphere in Section \ref{sec:spherical transport test}.
(Left) the initial condition, (centre) a computed state at $t=\tau/2$, as the hills have been deformed by the flow, and (right) the computed solution at $t=\tau$, as the tracers have returned to close to their initial condition.
Contours are spaced by $5\times10^{-3}$ kg kg$^{-1}$.
The superimposed arrows indicate the magnitude and direction of the transporting velocity field.}
\label{fig:aY_transport}
\end{figure}

\begin{figure}[h!]
\centering
\includegraphics[width=0.45\textwidth]{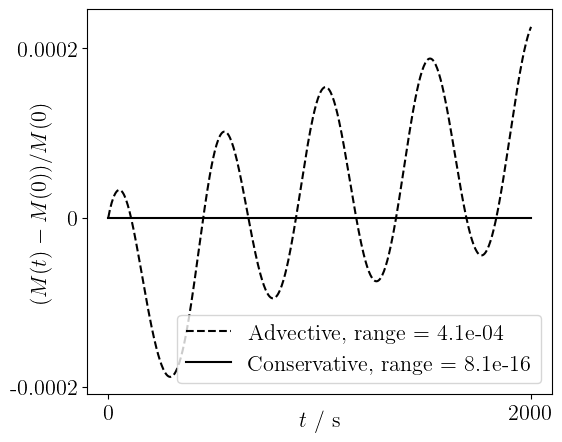}
\caption{Time series demonstrating the conservation of mass by tracer transport on a coarser mesh.
The evolution of tracer mass with initial conditions given by \eqref{eqn:conservation_config}, comparing the transport of a tracer with a purely advective transport scheme against the conservative transport described in Section \ref{sec:consistent transport}, showing that mass is indeed conserved in the latter case.
}
\label{fig:conservation_consistency}
\end{figure}

In this test case, a dry density field $\rho_d$ is transported by a prescribed wind on a fine resolution mesh, and a mixing ratio field $\overline{a}_Y$ is transported on a coarse mesh. This mimics the transport of tracers (for instance aerosols and chemicals) at a lower resolution, driven by a higher resolution dynamical core. This test is a variant of the time-dependent, deformational and divergent flow on the surface of the sphere from \citet{nair2010class}, \cite{lauritzen2012standard} to the sphere.
\\
\\
The spherical flow is defined as $u_0 = 2 \pi R/\tau$, $v_0=R/\tau$.
\begin{equation}
    u=u_0 \cos\left(\phi\right) - v_0 \cos\left(\frac{\pi t}{\tau}\right)\sin\left(\phi\right)\cos\left(\lambda - \frac{u_0 t}{R}\right)
\end{equation}
\begin{equation}
v=v_0 \cos\left(\frac{\pi t}{\tau}\right)\sin\left(\lambda - \frac{u_0 t}{R}\right),
\end{equation}
where $(\lambda,\phi)$ are the longitude and latitude, $R= 6.3781\times10^6 \;\mathrm{m}$ is the radius of the earth and $\tau=2000$ s is the length of the simulation. The finer and coarser meshes are C$32$ and C$16$ meshes respectively, where C$n$ denotes a cubed-sphere mesh with $n\times n$ cells per panel. In this case., the mesh is a two-dimensional spherical surface and the time step is $\Delta t = 4$ s.
\\
\\
The dry density initially varies with latitude while the tracer mixing ratio takes the form of two Gaussian hills. The initial conditions are
\begin{equation} \label{eqn:conservation_config}
\rho_d = \rho_0 + (\rho_t-\rho_0)\cos\phi, \qquad 
\overline{a}_Y = a_0 + a_t e^{-(L_{c_1}/R_{c_1})^2} + a_t e^{-(L_{c_2}/R_{c_2})^2},
\end{equation}
where $\rho_0=0.5$ kg m$^{-3}$, $\rho_t=1.0$ kg m$^{-3}$, $a_0=0.5$ kg kg$^{-1}$ and $a_t=1.0$ kg kg$^{-1}$.
$L_{c_1}$ and $L_{c_2}$ are the great circle distances between the local coordinate and the centre's of the bubbles $(\lambda_{c_1},\phi_{c_1})=(-\pi/4, 0)$, $(\lambda_{c_2},\phi_{c_2})=(\pi/4, 0)$ calculated as
\begin{equation}
L(\bm{x},\bm{x}_c)=\arccos[\sin(\phi)\sin(\phi_c)+\cos(\phi)\cos(\phi_c)\cos(\lambda - \lambda_c)],
\end{equation}
where $\bm{x}=(\lambda, \phi)$ and $\bm{x}_c=(\lambda_c, \phi_c)$.
The evolution of the mixing ratio $\overline{a}_Y$ is shown in Figure \ref{fig:aY_transport}.
This test can be used to demonstrate the conservation of mass of the tracer that is being transported on the coarse mesh.
Figure \ref{fig:conservation_consistency} shows time series of the tracer mass, comparing the approach described in \ref{sec:consistent transport} with an advective form of the transport equation, showing that mass is indeed conserved.
\\
\\

\subsection{Moist gravity wave} \label{sec:gravity wave test}
The next test is the moist gravity wave test case from \cite{bendall2020compatible}, adapted from the inertia-gravity wave test case of \cite{skamarock1994efficiency}. The final state in this test is spatially smooth, so it can be used to meaningfully measure the errors in the discretisation at different resolutions. This allows the effect of computing the physical parametrisation at a different resolution to be quantified. Like the rising bubble case of \cite{bryan2002benchmark}, the atmosphere is initially saturated and cloudy everywhere. Thus as air parcels move, water evaporates and condenses -- which is captured in our model through the use a Kessler physics scheme with no rain. This physics scheme can be performed on a different mesh to the dynamical core.
\\
\\
The domain is a two-dimensional vertical slice with height and length (10 km, 300 km). The initial conditions are the same as \cite{bendall2020compatible}, but with the exception that the definition of the wet equivalent potential temperature differs. In particular, LFRic-Atmosphere uses a latent heat $L_v$ which is constant with respect to temperature, and the heat capacities $c_p$ and $c_v$ used only the dry component of air. This means that the wet equivalent potential temperature is
\begin{equation} \label{equation:wet equivalent pot temp}
    \theta_e = \theta e^{L_v m_v / c_p T}.
\end{equation}
The same perturbation of and iterative procedure of \cite{bendall2020compatible} is used for the initial conditions. 
All runs used a mesh with $200$ vertical levels and a time step of $\Delta t= 1.2$ s.
\begin{figure}[h!]
\centering
\includegraphics[width=\textwidth]{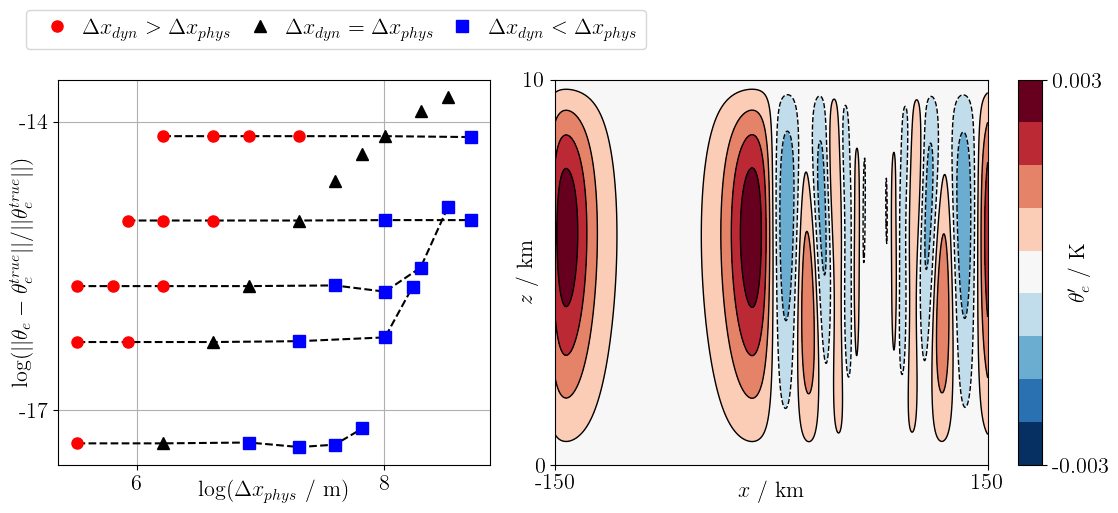}
\caption{(Left) Plots of the $L^2$ error norms in the final $\theta_e$ field from the moist gravity wave test of Section \ref{sec:gravity wave test}.
Errors are computed against a high resolution solution, and plotted as a function of the grid spacing of the physics mesh.
Dashed lines join the points corresponding to computations with the same resolution for the dynamical core.
The errors are largely independent of the resolution of the physics mesh, and instead depend strongly on the resolution of the dynamics mesh.
(Right) the final $\theta_e$ perturbation field, with both the dynamical core and the physics scheme using the same mesh. Contours are spaced at $3\times10^{-4}$ K.
}
\label{fig:moist_gw}
\end{figure}
The final state is shown in the right-hand side of Figure \ref{fig:moist_gw}. To investigate the effects using different meshes for the dynamical core and the physical parametrisation, we calculated the $L^2$ error norm of the final $\theta_e$ field relative to a high-resolution reference solution. The convergence plot in Figure \ref{fig:moist_gw} demonstrates that the errors are strongly dependent on the resolution of the dynamics mesh rather than the physics mesh. Increasing the horizontal resolution of the physics with respect to dynamics has no noticeable effect. Decreasing the resolution is shown to have a significant degradation in model solution quality only after a large enough resolution gap.

\subsection{Moist Baroclinic Wave} \label{sec:moist_baro_test}

The moist baroclinic wave test case forced by orography of \cite{hughes2023mountain} produces unstable waves and features important to the development of weather systems. Instead of a perturbation being added to the wind field, the orography forces the unstable atmospheric state and induces Rossby and inertia-gravity waves. The moist configuration of \cite{hughes2023mountain} is used, but with no-rain in the physics scheme. The test was run for $10$ days at C$96$, C$48$ and C$24$ resolutions with dynamics and physics at the same horizontal resolution. Two more configurations were run with dynamics at C$48$ but physics at C$96$ and C$24$. The grid is $30$ km in height with $30$ levels and a time step size of $\Delta t = 900$ s is used. The vertically stretched extrusion of \cite{ullrich2014proposed} is used.
\\
\begin{figure}[hp!]
\centering
\includegraphics[width=0.9\textwidth]{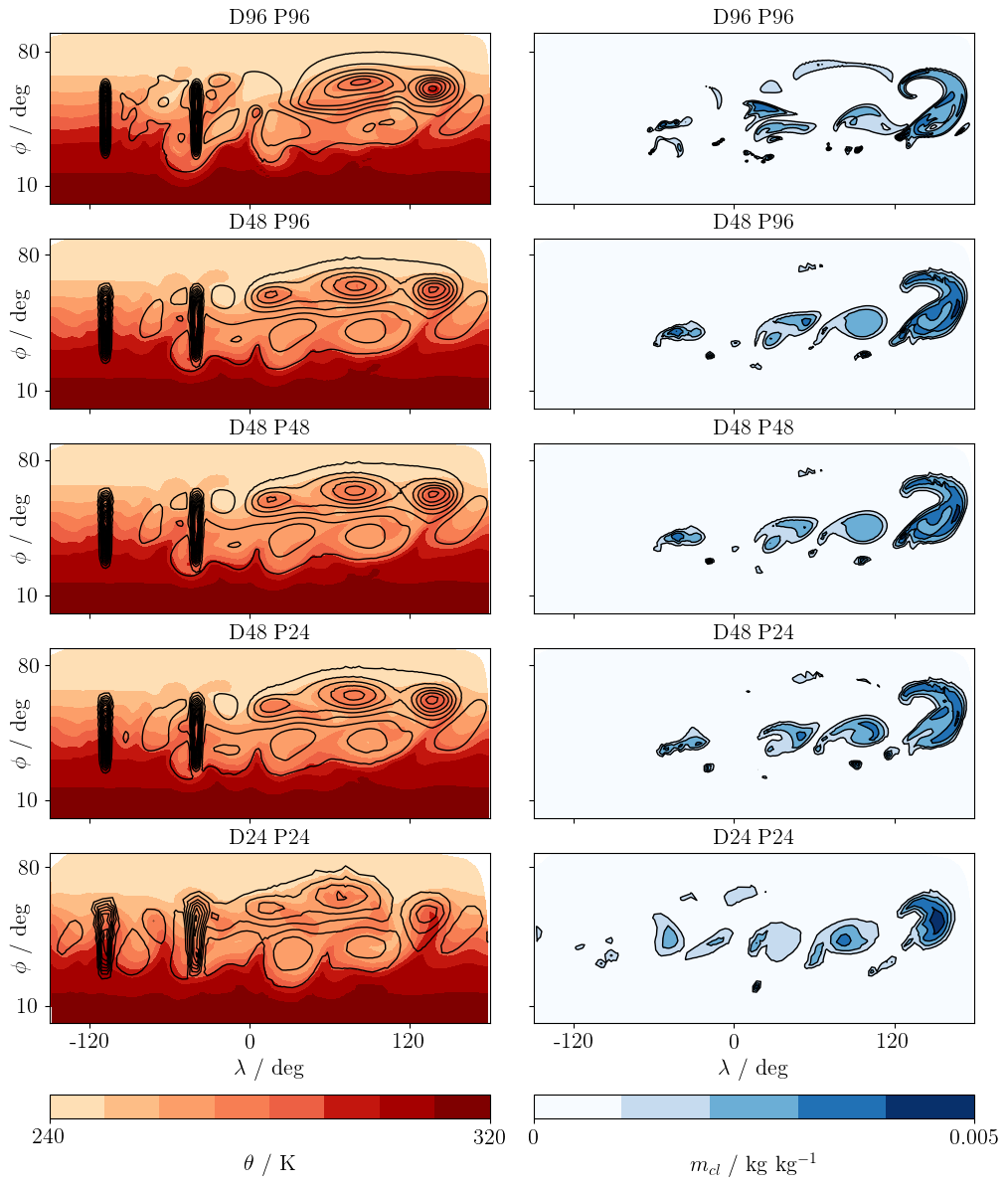}
\caption{Moist baroclinic wave test case, forced by orography at $t = 10$ days. The dynamics resolution is denoted by D$n$ and physics resolution is denoted P$n$.
(Left) Contours of the Exner pressure in the lowest model level, every $0.005$ (no unit) from $0.94$ to $1.01$, while the background colours show surface potential temperature with contours every 10 K.
(Right) Cloud liquid field at the $9$ km height with contours every $0.001$ kg kg$^{-1}$. In this test case the dynamical core resolution is much more influential in the evolution of the prognostic variables than the resolution of the physics scheme.}
\label{fig:moist_baro}
\end{figure} 
\\
The plots in Figure \ref{fig:moist_baro} display the first-level Exner and potential temperature fields, as well as the cloud at a height of $9$ km. Figure \ref{fig:moist_baro} shows that for this test case, increasing the physics resolution with respect to the dynamics has a much smaller impact than changing the resolution of the dynamical core. For the cases with the dynamical core using a C48 mesh, the pressure contours, temperature field and cloud fields look similar. Running the physics at a coarser resolution to the dynamics only degrades the solution minimally. It can be seen that when the dynamical core is run at the C$48$ resolution, there is a more prominent area of low pressure with strong pressure gradient at a latitude of around $140$ degrees; this drives cyclonic motion that results in the more overturned tail in the right-most cloud structure. This suggests that the cloud structure more strongly influenced by the fluid dynamics than the resolution of the latent heating effects.
\begin{figure}[htbp!]
\centering
\includegraphics[width=0.9\textwidth]{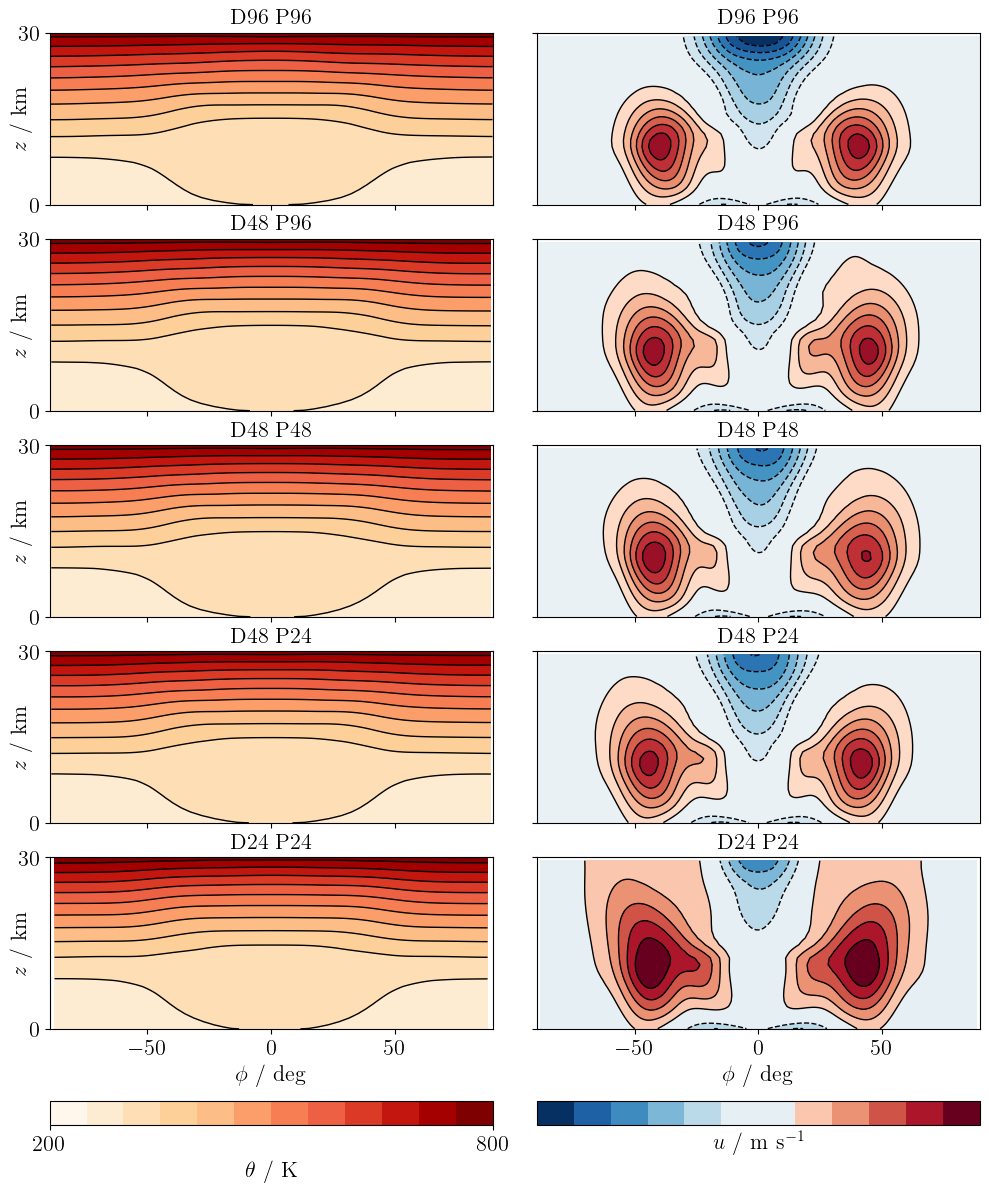}
\caption{The Held-Suarez test case with zonally averaged fields over $800$ days. The dynamics mesh is denoted by D$n$ and physics mesh is denoted by P$n$.
(Left) Zonally-averaged potential temperature, with contours every $50$ K, in a vertical slice.
(Right) Zonal wind with contours every $5$ m s$^{-1}$. The zero contour is omitted. The strength and extent of the jet appear to be largely dictated by the dynamics resolution, with D$48$ P$24$ appearing to be more comparable to D$96$ P$96$ and D$48$ P$48$ than D$24$ P$24$.}
\label{fig:held_suarez}
\end{figure}
\subsection{Held-Suarez} \label{sec:held_suarez}
The test case of \cite{held1994proposal} is a climate simulation which simulates the average atmospheric state by forcing the wind and surface temperature. This includes simple wind drag and temperature relaxation forcings, which can be treated as physical parametrisations, and so computed on a different mesh to the dynamical core. It generates two zonal jets in the mid-latitudes, and a vertical potential temperature gradient. For details of the set up of this test case in LFRic-Atmosphere, see \citet{sergeev2023simulations}. The test was run for $1000$ days, with a spin up time of $200$ days, meaning the results are averaged from the last $800$ days. The tests were run at differing time steps for different resolutions: 
C$96$ used a time step of $\Delta t = 900$ s, C$48$ had $\Delta t = 1800$ s and C$24$ used $\Delta t = 3600$ s. As with the baroclinic wave test, we also performed simulations with the dynamical core using a C$48$ mesh (and time step of $\Delta t = 1800$ s), but with the physical parametrisations on C$24$ and C$96$ meshes. The wind drag forcing is linear in the wind field, but is multiplied by a drag factor depending on the Exner field. We computed the drag factor on the physics mesh, and mapped this to the dynamical core mesh to multiply the wind field to get the resulting increment. A temporal off-centering of $\alpha = 0.55$ for the semi-implicit scheme was used. 
\\
\\
Typically with this test case the strength and extent of the jets are a result of the resolution. From Figure \ref{fig:held_suarez} we can see that the strength and extent of the D$48$ P$48$, D$48$ P$24$ and D$48$ P$96$ runs are comparable, implying that varying the resolution of the physics has minimal observable effect.

\section{Summary} \label{sec:summary}
This work has presented a formulation for mapping LFRic-Atmosphere's prognostic variables between meshes, to allow different components of the atmospheric model to use different meshes.
These meshes have the same vertical structure but different horizontal resolutions, with the resolution of the finer mesh such that its cells are nested within the cells of coarser meshes.
With this new capability, computational resources can be targeted towards the components that deliver the greatest impact on the model's accuracy.
At the same time, it may be possible to dramatically reduce the cost of some physical parametrisations without seeing a degradation in the quality of the solution.
The formulation is designed to possess a set of properties described in Section \ref{sec:desirable properties}, including mass conservation, preservation of constant mixing ratio fields and avoiding the generation of negative moisture concentrations.
The results in Section \ref{sec:results} demonstrate that the formulation does have these properties.
Tracers on a coarser mesh (representing chemicals and aerosols) are transported conservatively but such that constant mixing ratios are preserved.
Moisture species are mapped conservatively, without generating negative values and still preserving constant mixing ratios.
An idealised moist gravity wave test allowed quantification of the errors in the discretisation, which were largely independent of the resolution of the physics process.
\\
\\
The primary goals for future work are to apply this formulation to realistic NWP and climate models and to assess the scientific consequences of computing individual physical parametrisations at different resolutions to the dynamical core. The moist baroclinic wave and Held-Suarez test cases of Sections \ref{sec:moist_baro_test} and \ref{sec:held_suarez} present idealised versions of NWP and climate simulations; in these cases decreasing the physics resolution did not significantly degrade the solutions, but increasing the physics resolution offered no improvement in solution quality.
One particular target is for LFRic-Atmosphere to emulate the Junior-Senior capability of the UKESM model, in which the UKCA chemistry and aerosol component is performed on a coarser mesh.
Although it is possible to use a different mesh for each physical parametrisation, some schemes are more closely related and share auxiliary variables and so it may be appropriate for these schemes to share a mesh.
For instance, the radiation scheme interacts with the chemistry and aerosol variables, so we intend to explore using the same mesh for these components.
While the test cases in Section \ref{sec:results} did not reveal benefits from using a higher resolution mesh for the physical parametrisations, there may be clearer effects in ensemble simulations (e.g. with stochastic physics schemes), and more realistic configurations, particularly for interactions with the land surface through boundary-layer and convection processes.
One other interesting approach would use the same mesh for the physical parametrisations and dynamical core, but filter the prognostic fields that are passed to the physical parametrisations.
This would address the problem described by \citet{lander1997believable} of the errors at the smallest scales being amplified by the physical parametrisations.

\section*{Acknowledgements}
The work presented here was funded through Met Office work packages 2.2 and 3.3 of the ExCALIBUR research programme.
The authors would like to thank John Thuburn and Marc Stringer for useful conversations during this project, and Nigel Wood for his suggestions on improving the manuscript.
This work has also been facilitated by the many contributors to the LFRic-Atmosphere model and the underpinning LFRic-infrastructure, but particularly so by 
Ricky Wong for his design of the infrastructure for mapping wind fields between meshes.

\bibliography{main.bib}

\end{document}